\documentclass[hidelinks,onefignum,onetabnum]{siamart251216}

\usepackage{lipsum}
\usepackage{amsfonts}
\usepackage{graphicx}
\usepackage{epstopdf}
\usepackage{algorithmic}
\ifpdf
  \DeclareGraphicsExtensions{.eps,.pdf,.png,.jpg}
\else
  \DeclareGraphicsExtensions{.eps}
\fi

\usepackage{microtype}
\usepackage{subcaption}
\usepackage{booktabs}
\usepackage{hyperref}
\usepackage{amsfonts,mathrsfs}
\usepackage{verbatim}
\usepackage{acronym,color}
\usepackage{pifont}
\usepackage[normalem]{ulem}
\usepackage{yhmath}
\usepackage{amsmath}
\usepackage{amssymb}
\usepackage{mathtools}
\usepackage{booktabs}
\usepackage{array}
\usepackage{pifont}

\newsiamremark{remark}{Remark}
\newsiamremark{hypothesis}{Hypothesis}
\crefname{hypothesis}{Hypothesis}{Hypotheses}
\newsiamthm{claim}{Claim}

\newsiamremark{assumption}{Assumption}

\def\la{\langle}
\def\ra{\rangle}

\def\cF{{\cal F}}

\def\1{{\bf 1}}
\def\A{\mathscr{A}}

\def\ol{\bar}

\def\b{\beta}

\def\a{\alpha}

\def\o{\omega}

\def\R{\mathbb{R}}

\def\re{\mathbb{R}}

\def\dist{{\rm dist}}

\newcommand{\EXP}[1]{\mathsf{E}\!\left[#1\right] }

\def\argmin{\mathop{\rm argmin}}
\def\be{\begin{equation}}
\def\ee{\end{equation}}

\def\bit{\begin{itemize}}
\def\eit{\end{itemize}}

\newcommand{\cmark}{\ding{51}}
\newcommand{\xmark}{\ding{55}}
\newcolumntype{C}[1]{>{\centering\arraybackslash}m{#1}}
\newcommand{\remove}[1]{}

\headers{Randomized Feasibility-Update Algorithms for SVIs}{A. Chakraborty and A. Nedi\'c}

\title{Randomized Feasibility-Update Algorithms for Stochastic Variational Inequality Problems\thanks{Submitted to the editors on August 6, 2024. Revised on September 15, 2025 and September 20, 2026. This manuscript, in whole or in part, is not submitted, accepted, or published elsewhere, including conference proceedings.
\funding{This work was funded by the ONR award N00014-21-1-2242 and the NSF grant CIF-2134256.}}}

\author{Abhishek Chakraborty\thanks{School of Electrical Computer and Energy Engineering, Arizona State University, Tempe, Arizona, USA
  (\email{achakr61@asu.edu}, \email{Angelia.Nedich@asu.edu}).}
\and Angelia Nedi\'c\footnotemark[2]}

\usepackage{amsopn}

\ifpdf
\hypersetup{
  pdftitle={Randomized Feasibility-Update Algorithms for Stochastic Variational Inequality Problems},
  pdfauthor={A. Chakraborty and A. Nedi\'c}
}
\fi

\begin{document}

\maketitle

\begin{abstract}
    This paper considers stochastic monotone variational inequalities whose feasible region is the intersection of a (possibly infinite) number of convex functional level sets. A projection-based approach or direct Lagrangian-based techniques for such problems can be computationally expensive if not impossible to implement. To deal with the problem, we consider randomized methods that avoid the projection step on the whole constraint set by employing random feasibility updates. In particular, we propose and analyze modified stochastic Korpelevich and Popov methods for solving monotone stochastic VIs. We introduce a modified dual gap function and prove the convergence rates with respect to this function. We illustrate the performance of the methods in simulations on a zero-sum matrix game.
\end{abstract}


\begin{keywords}
Stochastic Variational Inequality, Randomized Feasibility Method, Modified Stochastic
Korpelevich Method, Modified Stochastic Popov Method
\end{keywords}

\begin{MSCcodes} 47N10, 49J40, 90C25, 68W20, 65K15
\end{MSCcodes}


\section{Introduction}\label{sec:Intro}

Variational inequalities (VIs) arise in a wide range of domains, such as economics~\cite{van2012dynamic}, optimization and machine learning~\cite{bertsekas2003convex, boyd2004convex, facchinei2003finite}, transportation science~\cite{acciaio2021cournot}, multi-agent games and reinforcement learning~\cite{lanctot2017unified}. VIs are characterized by a mapping
and an underlying set over which we seek to find a solution. Given a well-defined mapping and a feasible constraint set, VI problems can be solved with standard algorithms~\cite{solodov1999new, korpelevich1976extragradient, popov1980modification, tseng1995linear, malitsky2015projected} under suitable assumptions on the mapping. These algorithms typically employ projection onto the constraint set. This work considers a stochastic VI whose constraint set is specified as {\it the intersection of the level sets of many convex functions}, which precludes the direct use of projection. The motivation for considering the problems with infinitely many constraints such as stochastic noisy constraints of the form $\{x\in\re^{n}\mid g(x,\chi)\le 0 \}$ \cite{necoara2022stochastic}, where $\chi$ is a continuous random variable, comes from problems emerging in online settings dealing with multi-agent reinforcement learning with safety constraints \cite{gu2021multi, miryoosefi2019reinforcement}, or energy constraints~\cite{bai2024learning}, Nash-Cournot equilibrium problems in economics, constrained GANs~\cite{heim2019constrained}, and in almost-sure chance-constrained problems~\cite{fercoq2019almost}. To address such problems, this paper proposes methods based on random feasibility updates, which have been studied for optimization but remain less explored for VIs.

Related work in optimization with a large number of functional constraints dates back to~\cite{polyak1969minimization, polyak2001random}, which studied feasibility problems and developed randomized feasibility algorithms. Subsequent randomized methods for convex feasibility problems appear in~\cite{nedic2010random, necoara2021minibatch}. When each constraint is projection-friendly, one can project onto randomly sampled constraints at each iteration~\cite{liu2015averaging,wang2016stochastic}; such random projections have primarily been studied for finitely many constraints. Work~\cite{fercoq2019almost} addresses optimization with infinitely many constraints via duality and smoothing techniques, whereas~\cite{nedic2011random, nedic2019random, necoara2022stochastic, singh2024stochastic} use random feasibility updates to solve problems with infinitely, finitely, infinitely, and finitely many constraints, respectively. Other projection-free approaches for finitely constrained problems include Frank-Wolfe methods~\cite{woodstock2025splitting}, which rely on solving a linear minimization oracle (LMO), and ADMM~\cite{boyd2011distributed}, which requires solving subproblems to compute minimization. These extra computations, together with updating the dual multipliers of all the constraints per iteration, make these methods computationally cumbersome and impossible to implement for infinitely many constraints compared to the randomized feasibility methods.

\begin{table}[htbp]
\centering
\setlength{\tabcolsep}{3.5pt}
\renewcommand{\arraystretch}{1.35}
\footnotesize
\caption{Comparison of full projection-free methods for VIs under Lipschitz continuity and sometimes with an additional bounded variance term. The rates for
\cite{wang2015incremental} are best-iterate bounds to
stepsize-dependent neighborhoods, where $\epsilon$ denotes solution vicinity. For the other methods, rate is expressed using dual gap function, whereas for our case, we used a modified dual gap function.
Keys -- r-SPRG: Random-Stochastic Projected Reflected Gradient, r-SSE: Random-Stochastic
Subgradient Extragradient methods, ACVI: ADMM-based interior-point method for constrained VIs, Op. Ex.: Operator Extrapolation, P-D: Primal-Dual, Kor.: Korpelevich, Pv.: Popov, RF: Randomized Feasibility, M: Monotone, SM: Strongly Monotone, Det.: Deterministic, Stoc.: Stochastic, MDG: modified dual gap, Con.: Constraint, $\Pi_{X_a}[\cdot]$: Projection on individual sampled constraint set. Our methods use a few sampled constraints/iteration reducing the compute time compared to others.}
\label{tab:comparison}
\begin{tabular}{
    |C{0.95cm}
    |C{2.1cm}
    |C{1.3cm}
    |C{1.2cm}
    |C{0.9cm}
    |C{2.9cm}
    |C{1.8cm}|
}
\hline
\textbf{Works}
& \textbf{Method}
& \textbf{Mapping}
& \shortstack{\textbf{\# Con.}}
& $\boldsymbol{\Pi}_{X_a}[\cdot]$
& \textbf{Conv. Rate}
& \textbf{Con./itr.}
\\ \hline
\cite{wang2015incremental}
& Random Proj.
& \shortstack{SM, Stoc.}
& Finite
& \cmark
& \shortstack{
    $O(\alpha_k/\gamma_k) + \epsilon$
}
& Samples
\\ \hline
\cite{cui2021analysis}
& r-SPRG, r-SSE
& \shortstack{M, Stoc.}
& Finite
& \cmark
& \shortstack{
    $O(T^{-1/2})$
}
& Samples
\\ \hline
\cite{yang2022solving}
& ACVI
& \shortstack{M, Det.}
& Finite
& \xmark
& \shortstack{
    $O(T^{-1/2})$ (last itr.)\\
    $O(T^{-1})$ (average itr.)
}
& Full
\\ \hline
\cite{chavdarova2023primal}
& Inexact-ACVI
& \shortstack{M, Det.}
& Finite
& \xmark
& \shortstack{
    $O(T^{-1/2})$ (last itr.)
}
& Full
\\ \hline
\cite{boob2026first}
& Op. Ex. + P-D
& \shortstack{M, Stoc.}
& Finite
& \xmark
& \shortstack{
    $O(T^{-1/2})$
}
& Full
\\ \hline
\textbf{Ours}
& Kor./Pv. + RF
& \shortstack{M, Stoc.}
& Infinite
& \xmark
& \shortstack{
    $O(T^{-1/2})$ (inv. avg.)
}
& Samples
\\ \hline
\end{tabular}
\end{table}

In the context of VIs, a similar observation holds. Incremental and random projection algorithms have been studied~\cite{wang2015incremental} for strongly monotone Stochastic Variational Inequalities (SVIs). Subsequently, \cite{cui2021analysis} extended \cite{wang2015incremental} to monotone and weakly sharp SVIs. When projection-friendly sets are unavailable, primal-dual algorithms can be employed. Work~\cite{boroun2023projection} addresses nonconvex-concave saddle-point problems using a primal-dual conditional-gradient method, and \cite{alizadeh2024randomized} solves stochastic monotone Nash games via a randomized block primal-dual scheme. Both works consider settings with finitely many constraints. Additionally, \cite{boob2026first} studies primal-dual methods for deterministic and stochastic VIs with operator-extrapolation techniques. Beyond simple primal-dual methods, penalty-based approaches have been explored, for example, \cite{jordan2023first} for Nash equilibrium problems, and ADMM-based methods~\cite{chavdarova2023primal,yang2022solving} that use log-barrier functions for constrained VIs. These Lagrangian-based methods can become increasingly computationally demanding as the number of constraints grows. For the case where the constraint set is an infinite intersection induced by a stochastic function, $\cap_{\chi \in \mathcal{D}}\{x \mid g(x,\chi) \leq 0\}$, with $\chi$ as a random variable with support $\mathcal{D}$, work~\cite{tang2023optimal} solves the reformulated problem with constraint $\{ x \mid \EXP{\max(0,g(x,\chi))} \leq 0 \}$ by a stochastic proximal primal-dual method. The method requires inner- and outer-loops making it computationally demanding. The sets $\cap_{\chi \in \mathcal{D}}\{x \mid g(x,\chi) \leq 0\}$ and $\{ x \mid \EXP{\max(0,g(x,\chi))} \leq 0 \}$ coincide as long as $\mathcal{D}$ has a full support and the function $g(x,\cdot)$ is continuous for all $x$. However, \cite{tang2023optimal} contains a technical issue while assuming constraint qualification, which we describe in detail in Appendix~\ref{app_flaw}. Work~\cite{zhang2025primal} proposes a primal constrained-gradient method for monotone and strongly monotone VIs that identifies active constraints and solves a subproblem at each iteration, adding computational overhead. To our knowledge, none of these works develop methods for problems with a large (or infinite) number of constraints that are more computationally efficient than existing approaches. This current paper addresses these gaps for monotone SVIs. 

\noindent \textit{Contributions:}
(1) To the best of our knowledge, this is the first work to consider SVIs with potentially infinitely many constraints (Section~\ref{sec:problem}) and the first to study randomized methods for SVIs with functional constraints. We propose a modified stochastic Korpelevich (extragradient) method~\cite{korpelevich1976extragradient, juditsky2011solving} (Section~\ref{sec:Kor_algo}) and analyze it (Section~\ref{sec:Kor_analysis}). To further reduce operator evaluations per iteration, we also develop and study a modified stochastic Popov method~\cite{popov1980modification, chakraborty2025popov} (Section~\ref{sec:Popov_method}). Both methods handle (potentially infinitely many) functional constraints via random feasibility updates~\cite{polyak2001random, nedic2019random}. 
(2) We show, for the first time, \emph{per-iterate} convergence in expectation to the constraint set, with a geometric decay in the number of sampled constraints (Section~\ref{sec:rand_feas_updates}). By contrast, the randomized projection scheme of~\cite{cui2021analysis} achieves only an $O(1/\sqrt{k})$ rate, and the result of~\cite{wang2015incremental} is not per-iterate but for the minimum distance attained at an iteration, providing a rate only toward a ``vicinity of the set". 
(3) We handle monotone SVIs whose operator may exhibit discontinuous growth~\cite{juditsky2011solving, boob2026first}. We derive convergence rates for single-valued mappings (cf. Sections~\ref{sec:Kor_analysis} and~\ref{sec:Popov_method}) and show how the analysis extends to set-valued mappings satisfying an analogous condition (cf.\ Section~\ref{sec:multi_valued_F}). Unlike optimization, where objective values are available, we introduce a modified dual gap function for analyzing SVIs and prove an $O(1/\sqrt{k})$ rate under various iterate averaging schemes. The challenges in the analysis hinge on careful term decompositions and bounds on the infeasibility gap for averaged iterates under varying sample sizes in the randomized feasibility updates. Notably, both sample sizes and step sizes can be chosen to be problem-parameter-free (Remark~\ref{rem_Kor_parameter_free}, Theorem~\ref{thm_Popov}). Although primal-dual methods~\cite{boob2026first, chavdarova2023primal,yang2022solving} attain the same $O(1/\sqrt{k})$ rate, they are computationally heavier in large-constraint regimes because all constraints are processed in each iteration (cf.\  Section~\ref{sec:simulation}, Table~\ref{tab_runtime}).

\noindent
\textit{Notations:} We consider the vector space $\R^n$.  
 The inner product of two vectors $x$ and $y$ is denoted by $\la x,y\ra$, 
 while $\|\cdot\|$ is the standard Euclidean norm. The distance of a vector $\ol x$ 
 from a (nonempty) closed convex set $S$ is given by
 $\dist(\ol x,S)=\min_{x\in S}\|x-\ol x\|.$
 The projection of a vector $\ol x$
 on the set $S$ is  $\Pi_S[\ol x]=\argmin_{x\in S}\|x-\ol x\|^2.$
 For a scalar $a$, we use $a^+$ to denote the maximum of $a$ and~$0$,
 i.e.\ $a^+=\max\{a,0\}$.
 We use $\EXP{\o}$ for the expectation of a random variable~$\o$.
 We often abbreviate {\it almost surely} by {\it a.s.}


\section{Problem Formulation}\label{sec:problem}

Consider the VI problem of finding $x^* \in S := X \cap Y$
\begin{align}
    &\text{such that } \la F(x^*), x-x^* \ra \geq 0 \qquad \text{for all } x \in S = X \cap Y, \nonumber\\
    &\text{with } X := \cap_{a \in \A} X_{a} \quad \text{and } X_{a}:= \{x \in \re^{n}\mid g_{a}(x)\le 0\}, \label{VI_problem}
\end{align}
where $F: Y \rightarrow \R^n$ is the mapping of the VI. We view $\A$ as an index set (possibly infinite) and $a$ as its associated index element. One can also view the index $a$ as a random variable taking values in the set $\A$. In this case, given a random value $a$, 
the constraint set $X_{a}$ would coincide with $\{x \in \R^{n} \mid g(x,a) \leq 0\}$, similar to \cite{necoara2022stochastic}. It is assumed that the set $Y$ is simple to project on, but the projection on the set $S$ is computationally demanding, and we seek to reduce the computation by employing random projections. We present an assumption on the sets $X$ and $Y$.
\begin{assumption}\label{asum_closed_set}
The set $Y \subseteq \re^{n}$ is closed and convex, the set $S = X \cap Y$ is nonempty, and the function $g_{a}:\re^{n}\to\re$ is convex for all $a \in\A$. 
\end{assumption}
By Assumption~\ref{asum_closed_set}, each function $g_{a}$ is continuous over $\re^{n}$ \cite[Proposition 1.4.6]{bertsekas2003convex}, implying that 
 the set $\{x\in\re^{n}\mid g_{a}(x)\le0\}$ is closed and convex. Hence, the set $X$ is also closed and convex. Additionally, the subdifferential set $\partial g_{a}(x)$ is nonempty for all $x\in \re^{n}$ and all $a \in\A$~\cite[Proposition 4.2.1]{bertsekas2003convex}, implying that $\partial g_{a}^+(x)\ne\emptyset$
for all $x\in \re^{n}$ and $a \in\A$. The next assumption is on the boundedness of the set $Y$.
\begin{assumption}\label{asum_bounded_Y}
    The constraint set $Y$ is bounded, i.e., $\max_{x,y \in Y} \|x-y\|^2 \leq D$. 
    If the constraint index set $\A$ is infinite, assume that a scalar $M_g>0$ exists such that
  \begin{align}
      \|d\| \leq M_g \;, \quad  \forall d \in \partial g_{a}(x) \;, x \in Y \text{ and } a\in\A . \label{subgrad_norm_bd}
  \end{align}
\end{assumption}
When the set $Y$ is bounded and the constraint index set $\A$ is finite, the constraint function $g_{a} : \R^n \rightarrow \R$ has bounded subgradients on the set $Y$ \cite[Proposition 4.2.3]{bertsekas2003convex} for every $a\in\A$, implying that relation~\eqref{subgrad_norm_bd} is satisfied.
Assumption~\ref{asum_bounded_Y} ensures that the set $S$ in problem~\eqref{VI_problem} is bounded. The algorithms we study for \eqref{VI_problem} produce updates that are feasible for $Y$ but may be infeasible for $X$. We assume the mapping $F$ satisfies the following additional properties.
\begin{assumption}\label{asum_map_growth}
    The mapping $F:Y \rightarrow \R^{n}$ has the following variation over the set $Y$ with constants $L>0$ and $M \geq 0$,
   i.e.,
   \begin{align*}
       \|F(x)-F(y)\| \le L \|x-y\| + M \quad \forall x,y\in Y.
   \end{align*}
\end{assumption}
Assumption~\ref{asum_map_growth} allows discontinuities and applies to both single- and set-valued mappings; when $M=0$, $F$ is $L$-Lipschitz continuous. We consider single-valued mappings from here through Section~\ref{sec:Popov_method} and extend the framework to set-valued mappings in Section~\ref{sec:multi_valued_F}, with appropriate modifications to the notation, VI formulation, definitions, and assumptions, while retaining the same core ideas and analysis.
\begin{assumption}\label{asum_monotone}
    The mapping $F:Y \rightarrow \R^{n}$ is monotone
   over $Y$, i.e.,
   \begin{align}
       \la F(x) - F(y), x - y \ra \ge 0 \quad \forall x,y\in Y. \nonumber
   \end{align}
\end{assumption}

We conclude this section by presenting some preliminary results on sequences. 
\begin{lemma}\label{lem_seq}
    Let the sequence $\{\a_k\}$ be defined as $\a_k = \frac{\bar \alpha}{\sqrt{k+1}}$ for all $k \geq 0$, where $\bar \alpha > 0$ is a constant. Then for any $T \geq 1$, the following relations hold:
    \begin{align*}
        \sum_{k=1}^T \a_k \geq \frac{\bar \alpha \sqrt{T}}{\sqrt{2}} \quad \text{and} \quad \sum_{k=1}^T \a_k^2 \leq \bar \alpha^2 \ln(T+1) .
    \end{align*}
    Moreover, the sum of the reciprocal of the sequences is lower bounded as 
        \begin{align*}
            \sum_{k=1}^T \frac{1}{\a_k} \geq \frac{1}{\bar \alpha} \left[ \left( \frac{3}{2} \right)^{\frac{1}{2}} - \frac{2}{3} \right] T^{\frac{3}{2}} \quad \text{for any $T \geq 2$}.
        \end{align*}
\end{lemma}
\begin{proof}
    We see $\sum_{k=1}^T \a_k = \sum_{k=1}^T \frac{\bar \alpha}{\sqrt{k+1}} \geq \sum_{k=1}^T \frac{\bar \alpha}{\sqrt{T+1}} = \frac{\bar \alpha \sqrt{T}}{\sqrt{1+\frac{1}{T}}}$. Since $T \geq 1$, we obtain ${1+\frac{1}{T}} \leq 2$, which yields the first relation of the lemma. For the second relation, we upper estimate the summation using an integral, as follows
    $\sum_{k=1}^T \a_k^2 \leq \int_{k=0}^{T} \frac{\bar \alpha^2}{k+1} dk = \bar \alpha^2 \ln(T+1)$.

    To show the last relation of the lemma, we note that
    \begin{align}
        \sum_{k=1}^T \frac{1}{\a_k} \geq \int_{k=1}^T \frac{\sqrt{k+1}}{\bar \alpha} dk = \frac{2}{3 \bar \alpha} \left( (T+1)^{\frac{3}{2}} - 2^{\frac{3}{2}} \right) . \label{rel_lb}
    \end{align}
    Now, we want the following relation to hold for some constant $a > 0$, i.e.,
    \begin{align}
        (T+1)^{\frac{3}{2}} - 2^{\frac{3}{2}} \geq a T^{\frac{3}{2}} \quad \text{for any } T \geq 2 . \label{rel_lb2}
    \end{align}
    To do so, we define a function $\phi(T) = (T+1)^{\frac{3}{2}} - 2^{\frac{3}{2}} - a T^{\frac{3}{2}}$. We show that $\phi(T) > 0$ by proving that $\phi(T)$ is an increasing function of $T$ and it is non-negative for $T=2$. Taking the derivative of $\phi(T)$ with respect to $T$, we obtain $\frac{d\phi(T)}{d(T)} = \frac{3}{2} \left( (T+1)^{\frac{1}{2}} - a T^{\frac{1}{2}} \right) > 0$ for any $0< a \leq 1$ and $T \geq 2$. Letting $\phi(2) = 0$, we obtain $a = \left( \frac{3}{2} \right)^{\frac{3}{2}} - 1$. Hence, relation~\eqref{rel_lb2} holds, and the result follows by combining relations~\eqref{rel_lb} and~\eqref{rel_lb2}. 
\end{proof}


\section{Modified Stochastic Korpelevich Method}\label{sec:Kor_algo}

It is well known that the projection method may not converge for monotone mappings (Assumption~\ref{asum_monotone}) in general~\cite{mokhtari2020unified}. To circumvent this limitation, we employ the Korpelevich (a.k.a.\  Extragradient) method~\cite{korpelevich1976extragradient} for solving a monotone SVI problem with a mapping $F(\cdot) = \EXP{\widehat F(\cdot, \xi)}$, where $\xi$ is a random variable with a distribution on a set $\Xi$. As the exact evaluations of $F$ over $Y$ are computationally prohibitive, we rely on stochastic estimates $\widehat F(x, \xi)$ to design algorithms for solving SVI problem~\eqref{VI_problem}.

We consider a modification of the stochastic Korpelevich method in Algorithm~\ref{algo_Kor_method}. The algorithm is initialized at a randomly chosen point $x_0 \in Y$ with a distribution supported on the set $Y$. Under the boundedness of the set $Y$ (cf. Assumption~\ref{asum_bounded_Y}), $\EXP{\|x_0\|^2}$ is finite. At iteration $k \geq 1$, the method performs two updates with the same step size $\a_{k-1}$. Firstly, it updates the iterate $x_{k-1}$ using the mapping $\widehat F(x_{k-1}, \xi_k^1)$, with a random sample $\xi_k^1 \in \Xi$, to obtain an auxiliary point $u_{k}$. In the next step, the iterate $x_{k-1}$ is updated using the mapping $\widehat F(u_k, \xi_k^2)$, where $\xi_k^2 \in \Xi$ is the second random sample, to obtain the iterate $v_{k}$. Then, the randomized feasibility method (discussed later in Algorithm~\ref{algo_proj_steps}) is applied to reduce the infeasibility gap between $v_{k}$ and the set $X$, rendering  $x_{k}$ as the output at iteration $k$.
%
%
%
%
\begin{algorithm}
		\caption{Modified Stochastic Korpelevich Method}
		\label{algo_Kor_method}
		\begin{algorithmic}[1]
			\REQUIRE{ Initial iterate $x_{0}$, step size $\a_k$}
                \FOR{$k=1,\ldots$}
                \STATE \textbf{Sample} $\xi_k^1 \in \Xi$ and \textbf{update} $u_{k} = \Pi_{Y} \left[x_{k-1} - \a_{k-1} \widehat F(x_{k-1}, \xi_k^1)\right]$
                \STATE \textbf{Sample} $\xi_k^2 \in \Xi$ and \textbf{update} $v_{k} = \Pi_{Y} \left[x_{k-1} - \a_{k-1} \widehat F(u_{k}, \xi_k^2) \right]$
                \STATE \textbf{Call Randomized Feasibility Algorithm \ref{algo_proj_steps}:} Input $v_{k}$ and $N_{k}$; obtain $x_{k}$
                \ENDFOR
		\end{algorithmic}
\end{algorithm}	

For notational ease, we denote the stochastic errors, for all $k \geq 1$, as
\begin{align}
    b_k^1 = F(x_{k-1}) - \widehat F(x_{k-1}, \xi_{k}^1)  \qquad \text{and} \qquad b_k^2 = F(u_{k}) - \widehat F(u_{k}, \xi_{k}^2) . \label{errors} 
\end{align}

We next present a lemma for Algorithm~\ref{algo_Kor_method} showing a relation for the squared distance of the point $v_k$ from any arbitrary point $x \in Y$.

\begin{lemma}\label{lem_Kor1}
    Under Assumptions~\ref{asum_closed_set} and \ref{asum_map_growth}, and $\a_k > 0$, the following relation holds for the iterates $v_{k+1},x_k \in Y$ generated by Algorithm~\ref{algo_Kor_method}, for any point $x \in Y$, 
    \begin{align}
    \|v_{k+1} - x \|^2 &\leq \| x_{k} - x \|^2  + 2 \a_k \la b_{k+1}^2, u_{k+1} - x \ra + 2 \a_k \la F(u_{k+1}), x - u_{k+1} \ra \nonumber\\
    & + \frac{2 \a_k^2}{w_1} \left( \|b_{k+1}^2\|^2 + \|b_{k+1}^1\|^2 \right) + \frac{2 M^2 \a_k^2}{w_3} - \left( 1 - \frac{2 L^2 \a_k^2}{w_2} \right) \| x_{k} - u_{k+1} \|^2 \nonumber\\
    &- \left(1 - \frac{2w_1+w_2+w_3}{2} \right) \| v_{k+1} - u_{k+1} \|^2 , \nonumber
\end{align}
where $w_1,w_2$ and $w_3$ are positive constants, and $b_{k+1}^2$ and $b_{k+1}^1$ are defined in \eqref{errors}.
\end{lemma}
\begin{proof}
    Using the definition of the iterate $v_{k+1}$ and the non-expansiveness property of the projection operator, we obtain for any $x \in Y$,
    \begin{align}
    \|v_{k+1} - x \|^2 &\leq \| x_{k} - \alpha_k \widehat F(u_{k+1}, \xi_{k+1}^2) - x \|^2 - \| x_{k} - \alpha_k \widehat F(u_{k+1}, \xi_{k+1}^2) - v_{k+1} \|^2 \nonumber\\
    & = \| x_{k} - x \|^2 - \| v_{k+1} - x_{k} \|^2 + 2 \a_k \la \widehat F(u_{k+1}, \xi_{k+1}^2), x - v_{k+1} \ra . \label{Kor_main1}
    \end{align}
    The last quantity on the right-hand side of relation~\eqref{Kor_main1} can be written as
\begin{align}
    2 \a_k \la \widehat F(u_{k+1}, \xi_{k+1}^2), x - v_{k+1} \ra = 2 \a_k \la F(u_{k+1}), x - u_{k+1} \ra + 2 \a_k \la b_{k+1}^2, u_{k+1} - x \ra \nonumber\\
    + 2 \a_k \la \widehat F(x_{k}, \xi_{k+1}^1), u_{k+1} - v_{k+1} \ra + 2 \a_k \la \widehat F(u_{k+1}, \xi_{k+1}^2) - \widehat F(x_{k}, \xi_{k+1}^1), u_{k+1} - v_{k+1} \ra , \nonumber
\end{align}
where $b_{k+1}^2$ is defined in relation~\eqref{errors}. For the term
$\| v_{k+1} - x_{k} \|^2$ on the right-hand side of relation~\eqref{Kor_main1}, we have $\| v_{k+1} - x_{k} \|^2 = \| v_{k+1} - u_{k+1} \|^2 + \| x_{k} - u_{k+1} \|^2 - 2 \la x_{k} - u_{k+1}, v_{k+1} - u_{k+1} \ra$. Using these relations in \eqref{Kor_main1} yields
\begin{align}
    \|v_{k+1} &- x \|^2 \leq \| x_{k} - x \|^2 - \| v_{k+1} - u_{k+1} \|^2 - \| x_{k} - u_{k+1} \|^2 + 2 \a_k \la b_{k+1}^2, u_{k+1} - x \ra \nonumber\\
    & + 2 \a_k \la F(u_{k+1}), x - u_{k+1} \ra + 2 \la x_{k} - \a_k \widehat F(x_{k}, \xi_{k+1}^1) - u_{k+1}, v_{k+1} - u_{k+1} \ra \nonumber\\
    & + 2 \a_k \la \widehat F(x_{k}, \xi_{k+1}^1) - \widehat F(u_{k+1}, \xi_{k+1}^2) , v_{k+1} - u_{k+1} \ra. \label{Kor_main3}
\end{align}
By definition of $u_{k+1}$ in Algorithm~\ref{algo_Kor_method} and noting $v_{k+1} \in Y$, the sixth term on the right-hand side of relation~\eqref{Kor_main3} is non-positive due to the  projection properties, i.e.,
\begin{align}
    2 \la x_{k} - \a_k \widehat F(x_{k}, \xi_{k+1}^1) - u_{k+1}, v_{k+1} - u_{k+1} \ra \leq 0 . \label{eq_not_positive}
\end{align}
The last term on the right-hand side of relation~\eqref{Kor_main3} can be written as
\begin{align}
    2 \a_k &\la \widehat F(x_{k}, \xi_{k+1}^1) - \widehat F(u_{k+1}, \xi_{k+1}^2) , v_{k+1} - u_{k+1} \ra \nonumber\\
    & = 2 \a_k \la F(x_k) - F(u_{k+1}), v_{k+1} - u_{k+1} \ra + 2 \a_k \la b_{k+1}^2 - b_{k+1}^1, v_{k+1} - u_{k+1} \ra , \nonumber
\end{align}
where $b_{k+1}^1$ and $b_{k+1}^2$ are defined in relation~\eqref{errors}. The preceding relation can be upper bounded using Cauchy-Schwarz inequality, the growth condition on the mapping $F$ (Assumption~\ref{asum_map_growth}), and then Young's inequality, yielding the relation
\begin{align}
    &2 \a_k \la \widehat F(x_{k}, \xi_{k+1}^1) - \widehat F(u_{k+1}, \xi_{k+1}^2) , v_{k+1} - u_{k+1} \ra
    \leq \frac{2 \a_k^2}{w_1} \left( \|b_{k+1}^2\|^2 + \|b_{k+1}^1\|^2 \right) \nonumber\\
    & + \frac{2 L^2 \a_k^2}{w_2} \|x_k-u_{k+1}\|^2 + \frac{2 M^2 \a_k^2}{w_3} + \left(\frac{2w_1+w_2+w_3}{2} \right) \|v_{k+1} - u_{k+1}\|^2 , \nonumber
\end{align}
where $w_1, w_2$, and $w_3$ are some positive constants. Substituting the preceding relation back to relation~\eqref{Kor_main3} we obtain the stated relation.
\end{proof}

Lemma~\ref{lem_Kor1} is instrumental in establishing the convergence of Algorithm~\ref{algo_Kor_method}. To proceed with the proof, we must first derive a relation between the iterates $v_k$ and $x_k$ for all $k \geq 1$. This will be accomplished in the discussion of the randomized feasibility method and its associated analysis.


\section{Randomized Feasibility Updates}\label{sec:rand_feas_updates}

We consider the randomized feasibility updates following~\cite{polyak2001random,nedic2019random,necoara2021minibatch} to bypass the projection onto $X$, by randomly sampling sets $X_{a}$ from the family $\{X_{a}, a\in \A\}$ and performing sequential feasibility updates. Unlike~\cite{wang2015incremental, cui2021analysis}, which apply direct projection onto a randomly selected set, we use random feasibility updates since projections onto functional constraints may not have closed form expressions.

We next present a result for one step feasibility update for a functional constraint.
\begin{lemma}\label{lemma:basiter}
   Let $h$ be a convex function over a nonempty convex closed set $Z$. Given a vector $z\in Z$, a non-zero direction $d\in\partial h^+(z)$, and step size $\beta>0$, let $\hat z$ be given by $\hat z=\Pi_{Z}\left[z-\b \frac{h^+(z)}{\|d\|^2}\,d\right]$.
    Then, for any $\bar z\in Z$ such that $h^+(\bar z)=0$, we have
     \[\|\hat z -\bar z\|^2 \le  \|z -\bar z\|^2
       -\b(2-\beta)\,\frac{(h^+(z))^2}{\|d\|^2}.\]
\end{lemma}
The proof 
can be found in~\cite[Theorem~1]{polyak1969minimization}; also see~\cite{polyak2001random}. For $0<\beta<2$, Lemma~\ref{lemma:basiter}
implies that the point $\hat z$ is closer to the level set $\{\widetilde z \in Z \mid h(\widetilde z) \leq 0\}$ than the point $z$.
%
%

%
%
%
\begin{algorithm}
		\caption{Random Feasibility Steps}
		\label{algo_proj_steps}
		\begin{algorithmic}[1]
			\REQUIRE{$v_{k}$, $N_{k}$, deterministic step size $0 < \beta < 2$}
                \STATE \textbf{Set:} $z_{k}^0=v_{k}$
                \FOR{$i=1,\ldots,N_{k}$}
                \STATE \textbf{Sample:} Choose a random index $\o_{k}^{i} \in \A$
                \STATE \textbf{Compute:}
                Subgradient
                $d_{k}^{i}$ of  $g^+_{\o_{k}^{i}}(z)$ at the point $z=z_{k}^{i-1}$ 
                \STATE \textbf{Update:} $z_{k}^i
   = \Pi_{Y} \left[z_{k}^{i-1} - \beta\, \frac{g^+_{\o_{k}^{i}}(z_{k}^{i-1})}{\|d_{k}^{i}\|^2}\, d_{k}^{i}\right]$
                \ENDFOR
                \STATE \textbf{Output} $x_{k}=z_{k}^{N_{k}}$ .
		\end{algorithmic}
\end{algorithm}	
Algorithm~\ref{algo_proj_steps} performs feasibility updates for the functional constraints that define the set $X$. We combine this with variational inequality algorithms such as the Korpelevich method (Algorithm~\ref{algo_Kor_method}) and the Popov method, which will be introduced later in Algorithm~\ref{algo_Popov_method}. These algorithms generate updates using the mapping $F$. The inputs to Algorithm~\ref{algo_proj_steps} are the iterate $v_k$ (produced by the VI algorithms) and a deterministic time-varying batch size $N_k \geq 1$ at iteration $k$. The output of Algorithm~\ref{algo_proj_steps} is the updated iterate $x_k$.

With a deterministic step size $\beta\in(0,2)$, the Algorithm~\ref{algo_proj_steps} takes $N_{k}$ feasibility steps, each of which is random i.e.,
 the index variable $\o_{k}^{i} \in \A$ is random for all $i=1,\ldots, N_{k}$. The algorithm uses
 $d_{k}^{i}\in\partial g^+_{\o_{k}^{i}}(z_{k}^{i-1})$ if $g^+_{\o_{k}^{i}}(z_{k}^{i-1})>0$ and, otherwise, it can use $d_{k}^{i}=d$
 for some fixed vector $d \ne 0$ (if $g^+_{\o_{k}^{i}}(z_{k}^{i-1})=0)$. We note that the choice of the vector $d \ne 0$ is nonessential, since
 $z_{k}^{i}=z_{k}^{i-1}$ for any $d \ne 0$ due to
 $g^+_{\o_{k}^{i}}(z_{k}^{i-1})=0$. By Assumption \ref{asum_closed_set}, $\partial g_{\omega_{k}^i}^+(x)\ne\emptyset$ for all $x\in \re^{n}$, $\omega_{k}^i\in\A$, $k \geq 1$, and $i = 1, \ldots, N_{k}$. Therefore, the feasibility updates are well defined. Our next assumption deals with the random variables $\o_{k}^i$ similar to \cite{nedic2019random,necoara2022stochastic}.
 \begin{assumption}\label{asum-regularmod} 
There exists a constant $c>0$ such that, 
 \[\dist^2(x,S)\le c\,\EXP{(g^+_{\o_{k}^i}(x))^2}  \quad\hbox{for all } x \in Y,\quad i=1,\ldots, N_{k}, \quad k \geq 1. \]
 \end{assumption}
Regarding the random sampling of the indices $\o_{k}^i$, Assumption~\ref{asum-regularmod} allows for independent and identically distributed sampling, as well as sampling with a distribution conditional on the past samples drawn at iteration $k$. The constant $c$ in Assumption~\ref{asum-regularmod} exists, for example, if the index set $\A$ is bounded and there is a global error bound for the function $f(x)=\sup_{a\in\A} \{g_{a}(x) +\delta_{Y}(x)\}$, where $\delta_Y$ is the characteristic function of the set $Y$, taking value $0$ for all points inside $Y$ and $+\infty$ otherwise. A function $f(x)$ has a global error bound when there exists a constant $\bar \gamma>0$ such that $\dist(x,[f\le0])\le \bar \gamma f^+(x)$ for all $x$,
where $[f\le0]$ is the lower-level set of $f$ associated with the zero value. When $x \in Y$, then $\dist(x, S)\le \bar \gamma G^+(x)$, where $G^+(x) := \sup_{a\in\A} g_{a}^+(x)$~\cite[Section~2]{chakraborty2026randomized}. When the index set $\A$ is finite, the set $Y$ is polyhedral, and each $g_{a}$ is a linear function, then an error bound exists~\cite[Lemma 3]{hoffman2003approximate}, \cite[Theorem 3.1]{li2013global}. Also, such a bound exists for quadratic and more general convex functions under a Slater condition \cite{luo1994extension,lewis1998error,li2013global}, as well as for infinitely many linear constraints~\cite{hu1989approximate}. Moreover, \cite[Proposition 3]{bertsekas1999note} shows that an error bound holds if and only if the set of the optimal Lagrange multipliers are bounded. Work~\cite[Appendix~A]{chakraborty2026randomized} establishes this connection and shows that under strong Slater condition a global error
bound exists, and $\bar{\gamma}$ coincides with the optimal Lagrange
multiplier. For a finite index set $\A$, if a global error bound holds, then Assumption~\ref{asum-regularmod} holds with $c=\bar{\gamma}^2/(\min_{a\in\A}\mathbb{P}(\omega=a))$, where $c$ depends on the sampling distribution~\cite[Section 2.2]{nedic2011random}. Next, we present a lemma and a condition under which Assumption~\ref{asum-regularmod} holds for an infinite number of constraints.
\begin{lemma}\label{lem_inf_const}
    Suppose that the error bound $\dist(x,S)\leq\bar{\gamma} G^+(x)$ and the sampling condition $\mathbb{P}_{\omega\sim P} \left(g_{\omega}^{+}(x)\geq \rho G^+(x) \right) \geq p_0$ hold for all $x \in Y$, where the constants $\rho\in(0,1]$ and $p_0\in(0,1]$. Then, Assumption~\ref{asum-regularmod} holds with $c = \bar \gamma^2/(p_0 \rho^2)$.
\end{lemma}
\begin{proof}
    Define $E(x):=
    \left\{\omega\in \A: g_{\omega}^{+}(x)\geq\rho G^+(x) \right\}$ for all $x \in Y \setminus S$. We see that
    \begin{align*}
        \mathsf{E}_{\omega\sim P} \left[(g_{\omega}^{+}(x))^2 \right] \!\geq
        \mathsf{E}_{\omega\sim P}
        \left[(g_{\omega}^{+}(x))^2 \mathbf{1}_{E(x)}
        \right] \!\geq
        \rho^2 (G^+(x))^2
        \mathbb{P}_{\omega\sim P} (E(x) ) \geq
        p_0\rho^2 (G^+(x))^2.
    \end{align*}
    Combining this inequality with
    $\dist^2(x,S)\leq\bar{\gamma}^2(G^+(x))^2$ yields Assumption~\ref{asum-regularmod} with $c = \bar \gamma^2/(p_0 \rho^2)$. For $x\in S$, the result follows from $\dist(x,S)=0$.
\end{proof}
Next, we present two examples for an uncountable and a countable infinite constraint families, respectively, that satisfy the conditions of Lemma~\ref{lem_inf_const} and hence Assumption~\ref{asum-regularmod}.

\smallskip
\noindent\textbf{Example 1:} Let $Y=[-1,1]$, $\A=[0,1]$, $\omega\sim\operatorname{Unif}[0,1]$ and $g_a(x) = (1+a)x$ for all $a \in \A$. Then the feasible set $S = Y\cap\bigcap_{a\in\mathcal{A}}\{x:g_a(x)\leq 0\} = [-1,0]$. For $x > 0$, $G^+(x) = \sup_{a\in[0,1]}(1+a)x=2x$ and $\dist(x,S)=x=G^+(x)/2$. Hence, error bound holds for $\bar \gamma = 1/2$. Consider $\rho = 3/4$. We obtain $\mathbb{P}\left(g_{\omega}^{+}(x)\geq \rho G^+(x) \right) = \mathbb{P}\left((1+\omega)x\geq \frac{3}{2}x \right) = \mathbb{P}\left(\omega\geq\frac{1}{2}\right) = \frac{1}{2} = p_0$. Hence, by Lemma~\ref{lem_inf_const}, $c = 8/9$.

\smallskip
\noindent\textbf{Example 2:} Let $Y=[-1,1]$, $\A = \mathbb{N}$ (natural numbers), $\mathbb{P}(\omega=j)=2^{-j}$ for all $j \in \mathbb{N}$. Consider $g_j(x)=\left(1+{1}/{j}\right)x$ for all $j\in\mathbb{N}$. Here, the set $S = Y\cap\bigcap_{j\in\mathbb{N}} \{x:g_j(x)\leq0\} = [-1,0]$. Similar to Example~1, we can obtain $\bar \gamma = 1/2$. Considering $\rho = 3/4$ for all $x > 0$, we see $\mathbb{P}\left(g_\omega^+(x)\geq \rho G^+(x)\right) = \mathbb{P}\left(\left(1+\frac{1}{\omega}\right)x \geq \frac{3}{2}x \right) = \mathbb{P}(\omega\in\{1,2\}) = \frac{3}{4} = p_0$. Hence, we can obtain Assumption~\ref{asum-regularmod} with $c = 16/27$.



The works~\cite{yang2022solving, chavdarova2023primal} study ADMM based methods and do not prove the iterate convergence to the solution set. The work~\cite{boob2026first}, which analyzes the primal–dual method, shows an approximate feasibility rate of \(O(1/\sqrt{T})\) for \(\EXP{\|g^+(\bar x_T)\|}\), where \(\bar x_T\) is the averaged iterate and \(g^+(\bar x_T)\) concatenates \(g_a^+(\bar x_T)\) over all \(a\in\mathcal{A}\) (with \(\mathcal{A}\) finite). In our case, an error-bound assumption yields exponential decay of infeasibility, i.e., the expected distance from an iterate to the feasible set, as a function of the number of random samples \(N_k\) at each iteration \(k\ge 1\).


In the sequel, we define the sigma-algebra $\mathcal{F}_{k}$, for $k \geq 1$, with $\mathcal{F}_0 = \{x_0\}$, as:
  \begin{align} 
      \mathcal{F}_{k} = \mathcal{F}_0 \cup \{\omega_{t}^{i} \mid 1 \leq i \leq N_{t}, 1 \leq t \leq k\} \cup \{\xi_t^1 \cup \xi_t^2 \mid 1 \leq t \leq k+1 \} , \label{sigma_algebra}
  \end{align}
where $x_0$ denotes the (possibly random) initialization, $\omega_{t}^{i} \in \A$ are the random indices corresponding to sampled constraint sets from $X$, and $\xi_t^1, \xi_t^2 \in \Xi$ are the random variables governing the stochastic evaluations of the mapping $\widehat F(x_{t-1}, \xi_t^1)$ and $\widehat F(u_{t}, \xi_t^2)$, defined in Algorithm~\ref{algo_Kor_method}. Also, we will use the quantity $q\in(0,1]$, defined as follows:
    \begin{align}\label{quantity-q}
    q=\frac{\b(2-\beta)}{c M_g^2 }.
    \end{align}
We next present a basic relation for the iterates of Algorithm~\ref{algo_proj_steps}.
\begin{theorem}\label{thm_inf_updates}
    Under Assumptions \ref{asum_closed_set}, and \ref{asum_bounded_Y}, for 
   the iterates $x_{k}$ obtained via Algorithm \ref{algo_proj_steps} we have surely
    for all $x \in S = X \cap Y$ and $k \geq 1$,
    \[\|x_{k} - x\|^2 \le \|v_{k} - x\|^2
       - \frac{\b(2-\beta)}{M_g^2} \,\sum_{i=1}^{N_{k}} (g_{\o_{k}^{i}}^+(z_{k}^{i-1}))^2.\]
       Additionally, under Assumption~\ref{asum-regularmod}, we obtain almost surely the following for $k \geq 1$,
    \begin{align}
        &\frac{\b(2-\beta)}{M_g^2} \! \sum_{i=1}^{N_{k}} \EXP{(g_{\o_{k}^{i}}^+(z_{k}^{i-1}))^2 \!\mid \! \cF_{k-1}} \!\ge\! \left( (1\!-\!q)^{-{N_{k}}} - 1\right)  \EXP{\dist^2(x_{k},S) \!\mid \! \cF_{k-1}}, \nonumber\\
        &\EXP{\|x_{k}-x\|^2 \!\mid \! \cF_{k-1}} \! \leq \! \|v_{k}-x\|^2 \!-\! \left( (1-q)^{-N_{k}} \!-\! 1  \right) \EXP{\dist^2(x_{k},S) \mid  \cF_{k-1}} , \nonumber
    \end{align}
    where $0< \b < 2$, $M_g$ is the bound on subgradient norms (cf.\ relation~\eqref{subgrad_norm_bd}), and $q$ as given in relation~\eqref{quantity-q}, which satisfies $q\in(0,1)$.
\end{theorem}
\begin{proof}
We use the definition of $z_{k}^i$ in Algorithm \ref{algo_proj_steps} and Lemma~\ref{lemma:basiter}, with $Z=Y$. Thus, for all $x \in S = X \cap Y$, we obtain
for all  $i=1,\ldots,N_{k}$, 
    \begin{align}\label{eq_infes}
        \|z_{k}^i - x\|^2 
        &\le  \|z_{k}^{i-1} - x\|^2
       -\b(2-\beta)\,\frac{(g_{\o_{k}^{i}}^+(z_{k}^{i-1}))^2}{\|d_{k}^{i}\|^2}.
    \end{align}
    By summing the preceding relation over $i=1,\ldots,N_{k}$, and by using $z^0_{k}=v_{k}$, $z^{N_{k}}_{k}=x_{k}$, and $\|d_{k}^i\|^2 \leq M_g^2$
    for all $k$ (cf. relation~\eqref{subgrad_norm_bd}), we can obtain
\begin{align}
        \|x_{k} - x\|^2 \le \|v_{k} - x\|^2
       - \frac{\b(2-\beta)}{M_g^2} \,\sum_{i=1}^{N_{k}} (g_{\o_{k}^{i}}^+(z_{k}^{i-1}))^2, \label{x_bound}
    \end{align}
    which implies the first stated relation of the theorem.
Next, taking minimum with respect to $x \in S$ on both the sides of~\eqref{eq_infes}, along with $\|d_{k}^i\|^2 \leq M_g^2$, yields 
\begin{align*}
    \dist^2(z_{k}^i,S) \leq \dist^2(z_{k}^{i-1},S) - \frac{\b(2-\beta)}{M_g^2} (g_{\o_{k}^{i}}^+(z_{k}^{i-1}))^2 \quad \forall i=1,\ldots,N_{k}.
\end{align*}
Taking the conditional expectation in the preceding relation, given the past sigma-algebra $\cF_{k-1}$, we obtain almost surely for all $i=1,\ldots,N_{k}$,
\begin{align}
    \EXP{\dist^2(z_{k}^i,S) \mid \cF_{k-1}} \leq \EXP{\dist^2(z_{k}^{i-1},S) \mid \cF_{k-1}} - \frac{\b(2 - \beta )}{M_g^2} \EXP{(g_{\o_{k}^{i}}^+ (z_{k}^{i-1}))^2 \mid \cF_{k-1}} . \nonumber
\end{align}
By using Assumption~\ref{asum-regularmod} for the last term in the preceding relation, we obtain
     \begin{align}
         \EXP{(g_{\o_{k}^{i}}^+(z_{k}^{i-1}))^2 \mid \cF_{k-1}} &= \EXP{\EXP{(g_{\o_{k}^{i}}^+(z_{k}^{i-1}))^2 \mid \cF_{k-1},z_{k}^{i-1}}} \nonumber\\
         &\geq  \frac{1}{c} \EXP{\dist^2(z_{k}^{i-1},S) \mid \cF_{k-1}} \quad a.s. \text{ for all $i=1,\ldots,N_{k}$}. \label{eq-startexpest}
     \end{align}
Combining the preceding two relations, we see that almost surely for all $i=1, \ldots, N_{k}$,
\begin{align}
    \EXP{\dist^2(z_{k}^i,S) \mid \cF_{k-1}} \leq (1-q)\,
    \EXP{\dist^2(z_{k}^{i-1},S) \mid \cF_{k-1}}, 
\end{align} 
where $q$ is  given by~\eqref{quantity-q} and satisfies $q>0$ since $\b\in(0,2)$. 
Note that the preceding relation implies that we must have $1-q\ge0$. If $q=1$, it would follow that 
$ \EXP{\dist^2(z_{k}^i,S) \mid \cF_{k-1}} =0$ for all $i=1,\ldots,N_{k}$,
which would imply that $\dist^2(z_{k}^i,S)=0$ almost surely for all $i=1,\ldots,N_{k}$, which is highly unlikely case. Thus, without loss of generality, we assume that  $q<1$.
Using the definition of $x_{k}$, i.e., $x_{k}=z_{k}^{N_{k}}$,
      we have almost surely for all $i=1,\ldots,N_{k}$,
      \[\EXP{\dist^2(x_{k},S)\mid \cF_{k-1}}\le
      (1-q)^{N_{k}-i+1} \EXP{\dist^2(z_{k}^{i-1},S) \mid \cF_{k-1}} \qquad  {a.s.}\]
     Using~\eqref{eq-startexpest} in the preceding relation, we obtain
     for all $i=1,\ldots,N_{k}$,
     \[ \EXP{g_{\o_{k}^{i}}^+(z_{k}^{i-1}))^2 \mid \cF_{k-1}}
     \ge \frac{1}{c}\frac{1}{(1-q)^{N_{k}-i+1}}\EXP{\dist^2(x_{k},S) \mid \cF_{k-1}} \qquad  {a.s}.\]
     Summing over $i=1,\ldots,N_{k}$, and simplifying $\sum\limits_{i=1}^{N_{k}} \frac{1}{(1-q)^{N_{k}-i+1}} = \frac{1- (1-q)^{N_{k}}}{q(1-q)^{N_{k}}}$, we get
     \begin{align}
         \sum_{i=1}^{N_{k}} \EXP{g_{\o_{k}^{i}}^+(z_{k}^{i-1}))^2 \mid \cF_{k-1}}
     \ge \frac{1}{c}\frac{\left(1- (1-q)^{N_{k}}\right)}{q(1-q)^{N_{k}}}  \EXP{\dist^2(x_{k},S) \mid \cF_{k-1}} \;\; a.s. \nonumber
     \end{align}
     Since $q=\frac{\b(2-\beta)}{c M_g^2 }$,  we multiply by $\frac{\b(2-\beta)}{M_g^2 }$ both sides of the above inequality and simplify the relation to obtain almost surely the second relation of the theorem.

 The third relation of the theorem can be obtained by taking the conditional expectation on both sides of relation~\eqref{x_bound} with respect to $\cF_{k-1}$, and then using the second relation of the theorem to bound the last term on the right-hand side of relation~\eqref{x_bound}.
\end{proof}

Next, we present a lemma that provides a bound on the expected distance of the iterate $x_{k}$ from its projection on the set $S$.
\begin{lemma}\label{lem_inf_geom}
Under Assumptions \ref{asum_closed_set}, \ref{asum_bounded_Y}, and \ref{asum-regularmod}, and the constant $q$ defined in relation~\eqref{quantity-q}, for 
   the iterate $x_{k}$ obtained by Algorithm \ref{algo_proj_steps} we have almost surely,
    \begin{align*}
      &\EXP{\dist^2(x_{k},S) \mid \cF_{k-1}} \leq (1-q)^{N_{k}} {\|v_{k} - x\|^2} \quad \text{for all } x \in S= X \cap Y ,\\
      \text{and } &\EXP{\dist(x_{k},S) \mid \cF_{k-1}} \leq (1-q)^{\frac{N_{k}}{2}} {\|v_{k} - x\|} \quad \text{for all } x \in S= X \cap Y,
    \end{align*}
    where $q\in(0,1)$ as given in~\eqref{quantity-q}.
\end{lemma}
\begin{proof}
    Since $\dist^2(x_k,S) \leq \|x_k - x\|^2$ for any $x \in S $, from the third relation of Theorem~\ref{thm_inf_updates}, we obtain almost surely for all $x \in S$,
\begin{align*}
    \EXP{\dist^2(x_k,S) \mid \cF_{k-1}}  \le  \|v_{k} - x\|^2  - \left( (1-q)^{-{N_{k}}}-1\right)\EXP{\dist^2(x_{k},S) \mid \cF_{k-1}}, \nonumber
\end{align*}
    which yields the first result. Applying Jensen's inequality $\EXP{\dist(x_{k},S) \mid \cF_{k-1}} = \EXP{\sqrt{\dist^2(x_{k},S)} \mid \cF_{k-1}} \leq  \sqrt{\EXP{\dist^2(x_{k},S) \mid \cF_{k-1}}}$ to the first relation of the lemma, the second relation is obtained.
\end{proof}
Lemma~\ref{lem_inf_geom} will be useful while showing the convergence rates of the iterates to the solution set. Under Assumption~\ref{asum_bounded_Y}, $\EXP{\dist(x_{k},S) \mid \cF_{k-1}} \leq (1-q)^{\frac{N_{k}}{2}} \sqrt{D}$. Hence, the rate decays exponentially with respect to the number of samples $N_k$ for all $k \geq 1$.


\section{Convergence rate of Modified Stochastic Korpelevich Method}\label{sec:Kor_analysis}

In this section, we establish the convergence rate of Algorithm~\ref{algo_Kor_method}, building on Lemma~\ref{lem_Kor1} and the relation between the iterates $v_k$ and $x_k$ derived in Section~\ref{sec:rand_feas_updates}. We begin by presenting an auxiliary lemma that facilitates the subsequent analysis.
\begin{lemma}\label{lem_new_seq}
    Let $ \{h_k\}_{k=0}^T $, for all $T \geq 1$, be a sequence generated by the update
\[
h_{k+1} = \Pi_S \left[ h_k - \a_k b_{k+1}^2 \right], \quad \text{for } k = 0, 1, \ldots, T-1,
\]
with initial point $ h_0 \in S $. Then, for all $ k = 0, 1, \ldots, T-1 $,
\[
2 \a_k \left\langle b_{k+1}^2, h_k - x \right\rangle \leq \|h_k - x\|^2 - \|h_{k+1} - x\|^2 + \a_k^2 \|b_{k+1}^2 \|^2 \quad \text{for all } x \in S,
\]
where $ b_k^2 $ is defined as in relation~\eqref{errors}.
\end{lemma}
%
%
Lemma~\ref{lem_new_seq} is a special case of a proximal mapping with Bregman distances~\cite[Lemma 3 and Corollary 2]{juditsky2011solving}, \cite[Lemma 5]{chakraborty2025popov}.
The next lemma leverages Lemma~\ref{lem_new_seq} and the monotonicity of $F$ (cf. Assumption~\ref{asum_monotone}) to streamline Lemma~\ref{lem_Kor1}.
\begin{lemma}\label{lem_Kor2}
    Under Assumptions~\ref{asum_closed_set}, \ref{asum_bounded_Y}, \ref{asum_map_growth}, and \ref{asum_monotone}, and the scalars $\zeta>0$, $0<w_4<1$, $0< \b < 2$, $M_g$ defined in relation~\eqref{subgrad_norm_bd}, with the step size $0<\a_k \leq \frac{\sqrt{1-w_4}}{\sqrt{2}L}$, and the stochastic errors $b_{k+1}^2$ and $b_{k+1}^1$ defined in relation~\eqref{errors}, the iterates $v_{k+1} , x_k \in Y$ of Algorithm~\ref{algo_Kor_method} satisfy the following relation for all $x \in S$,
\begin{align*}
    2 \a_k &\la F(x), \Pi_S[x_k] - x \ra \leq \|v_{k} - x\|^2 - \|v_{k+1} - x \|^2 + \|h_k - x\|^2 - \|h_{k+1} - x\|^2 \\
    & - \frac{\b(2-\beta)}{M_g^2} \,\sum_{i=1}^{N_{k}} (g_{\o_{k}^{i}}^+(z_{k}^{i-1}))^2 + \zeta \dist^2(x_k, S) + 2 \a_k \la b_{k+1}^2, u_{k+1} - h_k \ra \\
    &+ \left[ \left( \frac{1}{\zeta} + \frac{1}{w_4} \right) \|F(x)\|^2 + 6M^2 + 7 \|b_{k+1}^2\|^2 + 6 \|b_{k+1}^1\|^2 \right] \a_k^2. 
\end{align*}
\end{lemma}
\begin{proof}
    The second term on the right-hand side of the relation in Lemma~\ref{lem_Kor1} can be decomposed as $2 \a_k \la b_{k+1}^2, u_{k+1} - x \ra = 2 \a_k \la b_{k+1}^2, h_k - x \ra + 2 \a_k \la b_{k+1}^2, u_{k+1} - h_k \ra$.
The first term on the right-hand side can be upper bounded using Lemma~\ref{lem_new_seq}, which when combined with Lemma~\ref{lem_Kor1} yields
\begin{align}
    \|v_{k+1} - &x \|^2 \leq \| x_{k} - x \|^2 + \|h_k - x\|^2 - \|h_{k+1} - x\|^2 + 2 \a_k \la b_{k+1}^2, u_{k+1} - h_k \ra \nonumber\\
    & + 2 \a_k \la F(u_{k+1}), x - u_{k+1} \ra + \left[ \frac{2M^2}{w_3} + \left( \frac{2}{w_1} + 1 \right) \|b_{k+1}^2\|^2 + \frac{2}{w_1} \|b_{k+1}^1\|^2 \right] \a_k^2 \nonumber\\
    & - \left( 1 - \frac{2 L^2 \a_k^2}{w_2} \right) \| x_{k} - u_{k+1} \|^2 - \left(1 - \frac{2w_1+w_2+w_3}{2} \right) \| v_{k+1} - u_{k+1} \|^2 . \label{Kor_main6}
\end{align}
The fifth term on the right-hand side can be upper bounded using the monotonicity of the mapping $F$ (Assumption~\ref{asum_monotone}); then by adding and subtracting $\Pi_S[x_k]$ and, then, applying Young's inequality with positive constants $\zeta$ and $w_4$ yields
\begin{align}
    2 \a_k \la F(u_{k+1}), x - u_{k+1} \ra &\leq 2 \a_k \la F(x), x - u_{k+1} \ra 
    \leq 2 \a_k \la F(x), x - \Pi_S[x_k] \ra \nonumber\\
    &+ \left[ \frac{1}{\zeta} + \frac{1}{w_4} \right] \|F(x)\|^2 \a_k^2 + \zeta \dist^2(x_k, S) + w_4 \|u_{k+1}-x_k\|^2 . \label{split_for_dg}
\end{align}
Substituting relation~\eqref{split_for_dg} in equation~\eqref{Kor_main6} and rearranging the terms yields
\begin{align}
    2 \a_k \la F(x), \Pi_S[x_k] \!-\! x \ra &\!\leq\! \| x_{k} \!-\! x \|^2 \!\!-\! \|v_{k+1} \!-\! x \|^2 \!+\! \|h_k \!-\! x\|^2 \!-\! \|h_{k+1} \!-\! x\|^2 \!+\! \zeta \dist^2(x_k, \!S) \nonumber\\
    &\hspace{-3.4cm}+\! 2 \a_k \la b_{k+1}^2, u_{k+1} \!-\! h_k \ra \!+\! \left[\! \left( \frac{1}{\zeta} \!+\! \frac{1}{w_4} \right)\! \|F(x)\|^2 \!+\! \frac{2M^2}{w_3} \!+\! \left( \frac{2}{w_1} \!+\! 1 \!\right) \!\|b_{k+1}^2\|^2 \!+\! \frac{2 \|b_{k+1}^1\|^2}{w_1}  \!\right] \!\a_k^2 \nonumber\\
    & \hspace{-2.5cm} - \left( 1 - w_4 - \frac{2 L^2 \a_k^2}{w_2} \right) \| x_{k} - u_{k+1} \|^2 - \left(1 - \frac{2w_1+w_2+w_3}{2} \right) \| v_{k+1} - u_{k+1} \|^2 . \label{Kor_main7}
\end{align}
To have $1 - w_4 - 2 L^2 \a_k^2/w_2\ge0$, the step size should satisfy
$\a_k \leq \frac{\sqrt{(1-w_4)w_2}}{\sqrt{2}L}$. We choose $w_2 = 1$, so that $0 < \a_k \leq \frac{\sqrt{1-w_4}}{\sqrt{2}L}$
with $0<w_4<1$.
We select $w_1 = w_3 = \frac{1}{3}$ so the last term on the right-hand side of relation~\eqref{Kor_main7} vanishes. Finally, using the first relation of Theorem~\ref{thm_inf_updates} to upper bound the first term on the right-hand side of~\eqref{Kor_main7}, upon grouping the terms, we obtain the desired relation.
\end{proof}
The relation of Lemma~\ref{lem_Kor2} includes the term $\|F(x)\|^2$ for $x \in S$, which can be bounded using a reference point $x_{ref} \in Y$, the boundedness of the set $Y$ (Assumption~\ref{asum_bounded_Y}), and the growth property of the mapping $F$ (Assumption~\ref{asum_map_growth}), as follows: for any $x\in Y$,
\begin{align}
    \|F(x)\| &\leq \|F(x) - F(x_{ref})\| + \|F(x_{ref})\| \leq B , \quad B = L \sqrt{D} + M + \|F(x_{ref})\|, \label{bounded_map_eq}
\end{align}
where $D$ is from Assumption~\ref{asum_bounded_Y}.
Thus, the same bound holds for any $x \in S = X \cap Y$. 
Next, we present an assumption for the stochasticity of the mapping $\hat F$.
\begin{assumption}\label{asum_stoc_error}
    The stochastic mapping evaluation $\widehat F(x, \xi)$ for all $x \in Y$ and $\xi \in \Xi$ is unbiased, and its variance is bounded by $\sigma^2$, i.e.,
    \begin{align*}
        &\EXP{\widehat F(x, \xi)} = F(x) , \quad \text{and} \quad \EXP{\|\widehat F(x, \xi) - F(x) \|^2} \leq \sigma^2 .
    \end{align*}
\end{assumption}
Assumption~\ref{asum_stoc_error} applies to both random draws of $\xi_k^{1}$ and $\xi_k^2$ in place of $\xi$ from the distribution $\Xi$ for all $k \geq 1$, in Algorithm~\ref{algo_Kor_method}.
%
%
%
%
%

Any point $x^* \in S$ that solves problem~\eqref{VI_problem} is referred to as a \emph{solution}. A \emph{Minty solution} is a point $x^* \in S$ such that $\la F(x), x - x^* \ra \geq 0$ for all $x \in S$.
Under the monotonicity of the mapping $F$ (Assumption~\ref{asum_monotone}), every solution is also a Minty solution~\cite[Lemma 2.2]{huang2023beyond}. Furthermore, when $M = 0$ in Assumption~\ref{asum_map_growth}, the mapping $F$ is continuous, in which case the two solution notions coincide~\cite[Lemma 1.5]{kinderlehrer2000introduction}.

We consider the dual gap function~\cite{marcotte1998weak, facchinei2003finite, nemirovski2004prox, huang2023beyond, chakraborty2024popov} defined as
\begin{align}
    G(y) = \max_{x \in S} \la F(x), y-x \ra \quad \text{for all } y \in S. \label{dg_gap}
\end{align}
The dual gap function satisfies $G(y) \geq 0$ for all $y \in S$, and $G(y) =0$ if and only if $y = x^*$ is a \textit{Minty} solution of the variational inequality~\cite{juditsky2016solving,chakraborty2024popov}. Since $S=X\cap Y$, the dual gap value $G(y)$ may be negative for some $y \in Y\setminus S$,  so it is not a valid merit function on the set $Y$.
For monotone mappings, prior work on similar problems using random projections~\cite{cui2021analysis}, primal–dual methods~\cite{boob2026first}, or ADMM-based methods~\cite{yang2022solving, chavdarova2023primal} have analyzed convergence via either the dual gap or the gap function, which can be negative when the point $y$ is infeasible.
To address this issue, we define an auxiliary quantity to track the progress of our algorithm, namely, the modified dual gap function, given by
\begin{align}
    \widehat G(y) = \left\lvert \max_{x \in S} \la F(x), y-x \ra \right\rvert \quad \text{for all } y \in Y. \label{mod_dual_gp}
\end{align}
The modified dual gap function $\widehat G(\cdot)$ is nonnegative by the definition. Moreover, if $y \in S$ is a Minty solution of $\text{VI}(S, F)$, then $\widehat G(y) = 0$. Note that if we analyzed our methods using the dual gap function (cf.\ \eqref{dg_gap}), we would need to provide both upper and lower bounds. By contrast, an upper bound on the modified dual gap function (cf.\ \eqref{mod_dual_gp}) suffices to capture the worst-case convergence rate for the method and the feasibility violation.
%
%
%
We analyze Algorithm~\ref{algo_Kor_method} via the averaged iterates
\begin{align}
    \widehat x_k = \frac{\sum_{t = 1}^k \gamma_t x_t}{\sum_{t = 1}^k \gamma_t}  \quad \text{and } \quad \widetilde x_k = \frac{\sum_{t = 1}^k \gamma_t \Pi_S[x_t]}{\sum_{t = 1}^k \gamma_t} , \label{avg_pt}
\end{align}
with some weights $\gamma_t > 0$ for all $1 \leq t \leq k$. 
Note that $\widehat x_k \in Y$ but need not lie in $S$, while $\widetilde x_k \in S$. Relating $\EXP{\widehat G(\widehat x_k)}$ to $\EXP{G(\widetilde x_k)}$ and using ${G(\widetilde x_k)} \geq 0$ yields
\begin{align}
    \EXP{\widehat G(\widehat x_k)} \leq \EXP{G(\widetilde x_k)} + \EXP{\text{V}(\widehat x_k)}, \; \text{with } \text{V}(\widehat x_k) = {\lvert \max_{x \in S} \la F(x), \widehat x_k - \widetilde x_k \ra \rvert} , \label{inf_rel}
\end{align}
where $\text{V}(\widehat x_k)$ denotes the infeasibility term associated with the iterate $\widehat x_k$, and thus depends on the weight sequence $\{\gamma_t\}_{t=1}^k$ and the deterministic sample sizes $\{N_t\}_{t=1}^k$. Next, we present a lemma to bound this quantity.
\begin{lemma}\label{lem_inf_gap_bd}
    Under Assumptions~\ref{asum_closed_set}, \ref{asum_bounded_Y}, \ref{asum-regularmod}, and \ref{asum_map_growth}, along with $q$, $B$ and $\mathrm{V}(\widehat x_k)$ defined in relations~\eqref{quantity-q}, \eqref{bounded_map_eq}, and \eqref{inf_rel} respectively, $D \geq \max_{x,y \in Y} \|x-y\|^2$, and $\Gamma(r) = \int_{0}^\infty p^{r-1} \exp(-p) dp$ being the gamma function, the following relations hold:\\
    i) For the sample size $N_k \geq 1$ for all $k \geq 1$, $\EXP{\mathrm{V}(\widehat x_k)} \leq B \sqrt{D} \left(\max_{1\leq t \leq k} (1-q)^{\frac{N_t}{2}} \right)$.
    ii) For $N_t = \lceil t^{\frac{1}{r}} \rceil$ for any $r \geq 1$, $\EXP{\mathrm{V}(\widehat x_k)} \leq B \sqrt{D} \left( \frac{\max_{1 \leq t \leq k} \gamma_t}{\sum_{t=1}^k \gamma_t} \frac{2^r \Gamma(r+1)}{\left(\ln \left(\frac{1}{1-q} \right)\right)^r} \right)$.
\end{lemma}
\begin{proof}
    {\it (i)} The expected infeasibility value $\EXP{\text{V}(\widehat x_k)}$ can be upper bounded using Cauchy–Schwarz inequality and relation~\eqref{bounded_map_eq}, yielding
    \begin{align}
        \EXP{\text{V}(\widehat x_k)} \leq \left(\max_{x \in S} \|F(x)\| \right) \EXP{\|\widehat x_k - \widetilde x_k\|} \leq B \frac{\sum_{t=1}^k \gamma_t \EXP{\dist(x_t,S)}}{\sum_{t=1}^k \gamma_t} ,  \label{terminfes_error}
    \end{align}
    where $\dist(x_t,S) = \|x_t - \Pi_S[x_t]\|$. The preceding relation can be further upper-bounded using the second relation of Lemma~\ref{lem_inf_geom} and Assumption~\ref{asum_bounded_Y}, yielding
    \begin{align}
        \EXP{\text{V}(\widehat x_k)} \leq B \sqrt{D} \frac{\sum_{t=1}^k \gamma_t (1-q)^{\frac{N_t}{2}}}{\sum_{t=1}^k \gamma_t} . \label{inf_up_bd1}
    \end{align} 
    Using $\sum_{t=1}^k \gamma_t (1-q)^{\frac{N_t}{2}} \leq \left( \sum_{t=1}^k \gamma_t \right) \left(\max_{1\leq t \leq k} (1-q)^{\frac{N_t}{2}} \right)$ in relation~\eqref{inf_up_bd1}, we obtain 
    the result.

    \noindent {\it (ii)} The numerator of relation~\eqref{inf_up_bd1} can be upper bounded in a different way, as follows: $\sum_{t=1}^k \gamma_t (1-q)^{\frac{N_t}{2}} \leq \left( \max_{1 \leq t \leq k} \gamma_t \right) \left( \sum_{t=1}^k (1-q)^{\frac{N_t}{2}} \right)$, thus yielding 
    \begin{align}
        \EXP{\text{V}(\widehat x_k)} \leq B \sqrt{D} \frac{\max_{1 \leq t \leq k} \gamma_t}{\sum_{t=1}^k \gamma_t} \left(\sum_{t=1}^k \left(\sqrt{1-q} \right)^{\lceil t^{\frac{1}{r}} \rceil}\right) , \label{inf_up_bd2}
    \end{align}
    where $N_t = \lceil t^{\frac{1}{r}} \rceil$ for any $r \geq 1$. Next, we upper bound $\sum_{t=1}^k \left(\sqrt{1-q} \right)^{\lceil t^{\frac{1}{r}} \rceil}$ as
    \begin{align*}
        &\sum_{t=1}^k \left(\sqrt{1-q} \right)^{\lceil t^{\frac{1}{r}} \rceil} \leq \sum_{t=1}^k \left(\sqrt{1-q} \right)^{t^{\frac{1}{r}}} \leq \int_{0}^{k} \left(\sqrt{1-q} \right)^{t^{\frac{1}{r}}} dt \\
        &= \int_{0}^{k} \exp \left( \ln{\left(\sqrt{1-q} \right)^{t^{\frac{1}{r}}}} \right) dt = \int_{0}^{k} \exp \left( - a t^{\frac{1}{r}} \right) dt, \;\; \text{where } a = \frac{1}{2} \ln \left( \frac{1}{1-q} \right).
    \end{align*}
    Using the change of variable $p = a t^{\frac{1}{r}}$, implying $t = \frac{p^r}{a^r}$ and $dt = \frac{r}{a^r} p^{r-1} dp$, we obtain
    \begin{align*}
        \sum_{t=1}^k \left(\sqrt{1-q} \right)^{t^{\frac{1}{r}}} \leq \frac{r}{a^r} \int_{0}^{ak^{\frac{1}{r}}} p^{r-1} \exp(-p) dp \leq \frac{r}{a^r} \int_{0}^\infty p^{r-1} \exp(-p) dp = \frac{r}{a^r} \Gamma(r) , 
    \end{align*}
    where $\Gamma(r)$ is the gamma function. By the properties of the gamma function, we have $r \Gamma(r) = \Gamma(r+1)$. Simplifying the preceding relation and using it back in relation~\eqref{inf_up_bd2}, the desired relation follows.
\end{proof}
Lemma~\ref{lem_inf_gap_bd}(i) shows geometric decay of the infeasibility term with respect to the worst-case sample size $\min_{1 \leq t \leq k} N_t$ encountered up to iteration $k$. This result is general and holds for any choice of $\{N_t\}_{t \geq 1}$. For the special schedule $N_t = \lceil t^{\frac{1}{r}} \rceil$ with $r \geq 1$, Lemma~\ref{lem_inf_gap_bd}(ii) yields a bound of order $O \left( \frac{\max_{1 \leq t \leq k} \gamma_t}{\sum_{t=1}^k \gamma_t} \right)$. With a suitable choice of weights $\gamma_t$ (specified later), we can obtain decaying rates in terms of the iteration index $k$. For specific choices of $r$ (e.g., $r=2,3$), the constants Lemma~\ref{lem_inf_gap_bd}(ii) can be significantly improved beyond the generic $\Gamma(r+1)$. However, since that quantity is just a constant that does not affect the convergence order, we are not deriving them here. Similar bounds hold for any schedule $N_k = f(k)$ that ensures a diminishing rate for $\EXP{\text{V}(\widehat x_k)}$. Moreover, for $N_k = \max(N, \lceil k^{\frac{1}{r}} \rceil)$ with fixed $N \geq 1$, we obtain the better of the two rates from Lemma~\ref{lem_inf_gap_bd}(i)–(ii), i.e.,
\begin{align*}
    \EXP{\mathrm{V}(\widehat x_k)} \leq B \sqrt{D} \min \left( (1-q)^{\frac{N}{2}}, \frac{\max_{1 \leq t \leq k} \gamma_t}{\sum_{t=1}^k \gamma_t} \frac{2^r \Gamma(r+1)}{\left(\ln \left(\frac{1}{1-q} \right)\right)^r} \right) .
\end{align*}



In the next two subsections, we present the convergence rate results for Algorithm~\ref{algo_Kor_method} under two different averaging schemes.




\subsection{Averaging of iterates with $\alpha_k$ as weights}

We consider the case when $\gamma_t = \a_t$ (cf.~\eqref{avg_pt}) for all $1 \leq t \leq k$, $k \geq 1$ and simplify the results of Lemma~\ref{lem_Kor2}.
\begin{lemma}\label{lem_Kor_normal_avg}
    Let Assumptions~\ref{asum_closed_set}, \ref{asum_bounded_Y}, \ref{asum_map_growth}, \ref{asum_monotone}, \ref{asum-regularmod}, and \ref{asum_stoc_error} hold. In Algorithm~\ref{algo_Kor_method}, let $\beta\in(0,2)$ and $0< \a_k \leq \frac{\sqrt{1-w_4}}{\sqrt{2}L}$, with $0<w_4<1$. Then, for all~$T \geq 1$,
    \begin{align*}
        \EXP{\widehat G(\widehat x_T)} \leq &\frac{D}{\sum_{k=1}^T \a_k} + \frac{\left[\left( \frac{1-q}{q} + \frac{1}{w_4} \right) B^2 + 6M^2 + 13 \sigma^2 \right] \sum_{k=1}^T \a_k^2}{2 \sum_{k=1}^T \a_k} + \EXP{\mathrm{V}(\widehat x_T)},
    \end{align*}
    where $\widehat x_T = \frac{\sum_{k=1}^T \a_k x_k}{\sum_{k=1}^T \a_k}$, the term $\mathrm{V}(\widehat x_T)$ is defined in relation~\eqref{inf_rel} (with $\gamma_k = \alpha_k$), while the constants $q$ and $B$ are defined in relations~\eqref{quantity-q} and~\eqref{bounded_map_eq}, respectively. 
\end{lemma}
\begin{proof}
    Summing the relation of Lemma~\ref{lem_Kor2}, over $k=1,\ldots,T$ yields for any  $T\ge 1$,
    \begin{align}
        2 \sum_{k=1}^T \a_k &\la F(x), \Pi_S[x_k] - x \ra \leq \|v_{1} - x\|^2 - \|v_{T+1} - x \|^2 + \|h_1 - x\|^2 - \|h_{T+1} - x\|^2 \nonumber\\
    & \hspace{-1cm} - \frac{\b(2-\beta)}{M_g^2} \, \sum_{k=1}^T \sum_{i=1}^{N_{k}} (g_{\o_{k}^{i}}^+(z_{k}^{i-1}))^2 + \zeta \sum_{k=1}^T \dist^2(x_k, S) + 2 \sum_{k=1}^T \a_k \la b_{k+1}^2, u_{k+1} - h_k \ra \nonumber\\
    &+ \sum_{k=1}^T \left[ \left( \frac{1}{\zeta} + \frac{1}{w_4} \right) \|F(x)\|^2 + 6M^2 + 7 \|b_{k+1}^2\|^2 + 6 \|b_{k+1}^1\|^2 \right] \a_k^2, \label{Kor_main8}
    \end{align}
    where $\zeta>0$ and $0<w_4<1$.
    The quantity on the left-hand side can be written as follows: $2 \sum_{k=1}^T \a_k \la F(x), \Pi_S[x_k] - x \ra = \left( 2 \sum_{k=1}^T \a_k \right) \la F(x), \widetilde x_T - x \ra$,
    where $\widetilde x_T = \frac{\sum_{k=1}^T \a_k \Pi_S[x_k]}{\sum_{k=1}^T \a_k}$. Combining the preceding relation with relation~\eqref{Kor_main8}, taking the maximum over $x \in S$, and applying Assumption~\ref{asum_bounded_Y} yields the following relation in terms of the dual gap function (cf.~\eqref{dg_gap})
    \begin{align}
        G(\widetilde x_T) \leq &\frac{D}{\sum_{k=1}^T \a_k} - \frac{\b(2-\beta)}{2 M_g^2 \sum_{k=1}^T \a_k} \, \sum_{k=1}^T \sum_{i=1}^{N_{k}} (g_{\o_{k}^{i}}^+(z_{k}^{i-1}))^2 \nonumber\\
        &+ \frac{\zeta}{2 \sum_{k=1}^T \a_k} \sum_{k=1}^T \dist^2(x_k, S) + \frac{1}{\sum_{k=1}^T \a_k} \sum_{k=1}^T \a_k \la b_{k+1}^2, u_{k+1} - h_k \ra \label{Kor_main10}\\
        & \hspace{-1cm} + \frac{1}{2 \sum_{k=1}^T \a_k} \sum_{k=1}^T \left[ \left( \frac{1}{\zeta} + \frac{1}{w_4} \right) \left(\max_{x \in S} \|F(x)\|^2 \right) + 6M^2 + 7 \|b_{k+1}^2\|^2 + 6 \|b_{k+1}^1\|^2 \right] \a_k^2 . \nonumber
    \end{align}
    Next, we take the total expectation on both sides. Applying relation~\eqref{bounded_map_eq} and Assumption~\ref{asum_stoc_error}, the last term on the right-hand side can be upper bounded. The expected value of the second quantity on the right-hand side can be simplified using law of iterated expectation and Assumption~\ref{asum_stoc_error} as follows $\EXP{\sum_{k=1}^T \a_k \la b_{k+1}^2, u_{k+1} - h_k \ra} = \EXP{\sum_{k=1}^T \a_k \la \EXP{b_{k+1}^2 \mid \cF_{k-1} \cup \{\xi_{k+1}^1\}}, u_{k+1} - h_k \ra} = 0$. The second quantity on the right-hand side of~\eqref{Kor_main10} can be simplified using the law of iterated expectation, and by applying the second relation of Theorem~\ref{thm_inf_updates}, we obtain 
    \begin{align}
        - \EXP{\sum_{k=1}^T  \frac{\b(2-\beta)}{ M_g^2} \sum_{i=1}^{N_{k}} (g_{\o_{k}^{i}}^+(z_{k}^{i-1}))^2}
        \leq -\sum_{k=1}^T \left( (1-q)^{-{N_{k}}} - 1\right) \EXP{\dist^2(x_{k},S) }, \nonumber
    \end{align}
    where the constant $q$ is defined in relation~\eqref{quantity-q}. Substituting all the preceding relations back in relation~\eqref{Kor_main10} and grouping  the terms accordingly, we obtain
    \begin{align}
        \EXP{G(\widetilde x_T)} \leq \frac{D}{\sum_{k=1}^T \a_k} &- \frac{1}{2 \sum_{k=1}^T \a_k} \sum_{k=1}^T \left( (1-q)^{-{N_{k}}} - (1 + \zeta) \right) \EXP{\dist^2(x_{k},S)} \nonumber\\
        & + \frac{\left( \frac{1}{\zeta} + \frac{1}{w_4} \right) B^2 + 6M^2 + 13 \sigma^2}{2 \sum_{k=1}^T \a_k} \sum_{k=1}^T \a_k^2 . \label{Kor_main12}
    \end{align}
    We require $(1-q)^{-{N_{k}}} - (1 + \zeta) \geq 0$ for all $k \geq 1$, which yields $N_k \geq \frac{\log(1 + \zeta)}{\log\left( \frac{1}{1-q} \right)}$ for all $k \geq 1$. For the algorithm to work with any number of samples $N_k \geq 1$, we get $1 + \zeta = \frac{1}{1-q}$, which yields the value of the constant $\zeta= \frac{q}{1-q}$.
    Substituting the value of $\zeta$ back in relation~\eqref{Kor_main12} and dropping the second term to its right, we obtain 
    \begin{align}
        \EXP{G(\widetilde x_T)} \leq \frac{D}{\sum_{k=1}^T \a_k} + \frac{\left( \frac{1-q}{q} + \frac{1}{w_4} \right) B^2 + 6M^2 + 13 \sigma^2}{2 \sum_{k=1}^T \a_k} \sum_{k=1}^T \a_k^2 . \nonumber
    \end{align}
    Applying relation~\eqref{inf_rel} ($\gamma_k = \alpha_k$) to the preceding relation yields the final relation. 
\end{proof}
Lemma~\ref{lem_Kor_normal_avg} establishes a basic convergence result that depends on the choice of step size and certain associated constants. Furthermore, the result includes an infeasibility term $\EXP{\mathrm{V}(\widehat x_T)}$ that depends on the sample size selection. 
Next, we present a theorem based on a decaying step size scheme along with the expected infeasibility gap $\EXP{\text{V}(\widehat x_k)}$.
\begin{theorem}\label{thm_Kor1}
    Let Assumptions~\ref{asum_closed_set}, \ref{asum_bounded_Y}, \ref{asum_map_growth}, \ref{asum_monotone}, \ref{asum-regularmod}, and \ref{asum_stoc_error} hold, along with the step sizes $0<\beta<2$ and $0< \a_k = \min \left( \frac{\bar \alpha}{\sqrt{k+1}},\frac{\sqrt{1-w_4}}{\sqrt{2}L} \right)$ for $k \geq 1$, with constants $0<w_4<1$ and $\bar \alpha>0$. Then, for the iterates of Algorithm~\ref{algo_Kor_method}, we have for $T \geq 1$,
    \begin{align*}
        \EXP{\widehat G(\widehat x_T)} \leq &\frac{\sqrt{2} D}{\bar \alpha \sqrt{T}} + \frac{\frac{\bar \alpha}{\sqrt{2}} \left[\left( \frac{1-q}{q} + \frac{1}{w_4} \right) B^2 + 6M^2 + 13 \sigma^2 \right] \ln(T+1)}{\sqrt{T}} + \EXP{\text{V}(\widehat x_k)}, 
    \end{align*}
    where $\widehat x_T = \frac{\sum_{k=1}^T \a_k x_k}{\sum_{k=1}^T \a_k}$, $D \geq \max_{x,y \in Y} \|x-y\|^2$, $\mathrm{V}(\widehat x_T)$, $B$, and $q$ are defined in relations~\eqref{inf_rel}, \eqref{bounded_map_eq}, and \eqref{quantity-q}, respectively. Further, if $0< \a_k = \min \left( \frac{\bar \alpha}{\sqrt{T}},\frac{\sqrt{1-w_4}}{\sqrt{2}L} \right)$, then
    \begin{align*}
        \EXP{\widehat G(\widehat x_T)} \leq &\frac{D}{\bar \alpha \sqrt{T}} + \frac{\bar \alpha \left[\left( \frac{1-q}{q} + \frac{1}{w_4} \right) B^2 + 6M^2 + 13 \sigma^2 \right]}{\sqrt{T}} + \EXP{\text{V}(\widehat x_k)}. 
    \end{align*}
\end{theorem}
\begin{proof}
    If the step size satisfies $\alpha_k = \min \left( \frac{\bar \alpha}{\sqrt{k+1}}, \frac{\sqrt{1 - w_4}}{\sqrt{2}L} \right)$, then the relation in Lemma~\ref{lem_Kor_normal_avg} can be upper bounded using Lemma~\ref{lem_seq}. To obtain the second relation, we choose the constant $\alpha_k = \min \left( \frac{\bar \alpha}{\sqrt{T}}, \frac{\sqrt{1 - w_4}}{\sqrt{2}L} \right)$ and apply this choice to the relation in Lemma~\ref{lem_Kor_normal_avg}, which yields the desired result.
\end{proof}
Theorem~\ref{thm_Kor1} establishes an $O(\ln(T)/\sqrt{T})$ rate with diminishing step sizes; with a constant step size, the $\ln(T)$ factor is eliminated. An additional infeasibility term $\EXP{\mathrm{V}(\widehat x_T)}$ appears and we analyze its convergence rate in the next remark.
%
%
%
%
%
\begin{remark}\label{rem_sample_root_growth}
    For any sample size sequence $\{N_k\}_{k=1}^T$, by Lemma~\ref{lem_inf_gap_bd}(i) we obtain $\EXP{\text{V}(\widehat x_T)} \leq B \sqrt{D} \left(\max_{1 \leq k \leq T} (1-q)^{\frac{N_k}{2}}  \right)$, i.e., the term decays geometrically fast with respect to $\min_{1 \leq k \leq T} N_k$. When $N_k = \lceil k^{\frac{1}{r}} \rceil$ with any $r \geq 1$, then with $\gamma_k = \alpha_k = \min \left( \frac{\bar \alpha}{\sqrt{k+1}}, \frac{\sqrt{1 - w_4}}{\sqrt{2}L} \right)$ in Lemma~\ref{lem_inf_gap_bd}(ii) yields for any $0< w_4 < 1$, $\bar \alpha > 0$, and $T \geq 1$,
    \begin{align*}
        \EXP{\mathrm{V}(\widehat x_T)} \leq B \sqrt{D} \left( \frac{\min \left( \frac{\bar \alpha}{\sqrt{2}}, \frac{\sqrt{1 - w_4}}{\sqrt{2}L} \right)}{T \min \left( \frac{\bar \alpha}{\sqrt{T+1}}, \frac{\sqrt{1 - w_4}}{\sqrt{2}L} \right)} \frac{2^r \Gamma(r+1)}{\left(\ln \left(\frac{1}{1-q} \right)\right)^r} \!\right) .
    \end{align*}
    For $T \geq 1$, the quantity $\frac{\bar \alpha T}{\sqrt{T+1}} = \sqrt{T} \frac{\bar \alpha}{\sqrt{1+\frac{1}{T}}} \geq \frac{\bar \alpha \sqrt{T}}{\sqrt{2}}$. Hence, we obtain
    \begin{align*}
        \EXP{\mathrm{V}(\widehat x_T)} \leq \begin{cases}
            B \sqrt{D} \left( \frac{\min \left( \frac{\bar \alpha L}{\sqrt{1-w_4}}, 1 \right)}{T} \frac{2^r \Gamma(r+1)}{\left(\ln \left(\frac{1}{1-q} \right)\right)^r} \right) \quad &\text{if $\frac{\bar \alpha \sqrt{T}}{\sqrt{2}} \geq \frac{\sqrt{1 - w_4} T}{\sqrt{2}L}$,} \\
            B \sqrt{D} \left( \frac{\min \left( 1, \frac{\sqrt{1 - w_4}}{\bar \alpha L} \right)}{\sqrt{T}} \frac{2^r \Gamma(r+1)}{\left(\ln \left(\frac{1}{1-q} \right)\right)^r} \right) \quad &\text{if $\frac{\bar \alpha \sqrt{T}}{\sqrt{2}} < \frac{\sqrt{1 - w_4} T}{\sqrt{2}L}$.}
        \end{cases}
    \end{align*}
     Hence, a worst case convergence rate of $O \left( \frac{1}{\sqrt{T}} \right)$ is always achieved for the quantity $\EXP{\widehat G(\widehat x_T)}$. Moreover, for the constant step size $\a_k = \min \left( \frac{\bar \alpha}{\sqrt{T}},\frac{\sqrt{1-w_4}}{\sqrt{2}L} \right)$, the same analysis can be carried out and we can obtain $\EXP{\mathrm{V}(\widehat x_T)} \leq \frac{B \sqrt{D}}{T} \frac{2^r \Gamma(r+1)}{\left(\ln \left(\frac{1}{1-q} \right)\right)^r}$. Hence, a smaller constant step size $\a_k$ is beneficial for a faster infeasibility decay.
\end{remark}
Theorem~\ref{thm_Kor1} and Remark~\ref{rem_sample_root_growth} show an $O(1/\sqrt{T})$ rate can be achieved.
The dependence of $\ln(T)$ in Theorem~\ref{thm_Kor1} can be removed by choosing $\a_k = O(\bar \alpha /\sqrt{T})$, but such small steps may be undesirable for large $T$. Alternatively, the $\ln(T)$ factor can be eliminated by averaging the iterates with weights proportional to $\alpha_k^{-1}$, which we present next.



\subsection{Averaging of iterates with $\alpha_k^{-1}$ as weights}
We present a theorem that shows the convergence rate for the case when the weight $\gamma_k = \a_k^{-1}$ for all $k \geq 1$.
\begin{theorem}\label{thm_Kor2}
    Under Assumptions~\ref{asum_closed_set}, \ref{asum_bounded_Y}, \ref{asum_map_growth}, \ref{asum_monotone}, \ref{asum-regularmod}, and \ref{asum_stoc_error}, and with step size selections $0<\beta<2$ and $0< \a_k = \min \left( \frac{\bar \alpha}{\sqrt{k+1}}, \frac{\sqrt{1 - w_4}}{\sqrt{2}L} \right)$, $k \geq 1$, where the constants $\bar \alpha>0$ and $0<w_4<1$, then the iterates of Algorithm~\ref{algo_Kor_method} satisfy for $T \geq 2$,
    \begin{align*}
        &\EXP{\widehat G(\widehat x_T)} \leq \frac{\left( \frac{3 D}{\bar \alpha} + \bar \alpha \left[ \left( \frac{1-q}{q} + \frac{1}{w_4} \right) B^2 + 6M^2 + 13 \sigma^2 \right] \right)}{2 \left[ \left( \frac{3}{2} \right)^{\frac{1}{2}} - \frac{2}{3} \right] \sqrt{T}} + \EXP{\mathrm{V}(\widehat x_T)}
    \end{align*}
    where $\widehat x_T = \frac{\sum_{k=1}^T \a_k^{-1} x_k}{\sum_{k=1}^T \a_k^{-1}}$, the quantities $q$, $B$, $\mathrm{V}(\widehat x_T)$, and $D$ are defined in relations~\eqref{quantity-q}, \eqref{bounded_map_eq}, \eqref{inf_rel} and Assumption~\ref{asum_bounded_Y}, respectively.
\end{theorem}
\begin{proof}
    We divide both sides of the relation of Lemma~\ref{lem_Kor2} by $\a_k^{2}$, and add and subtract the quantity $\a_{k-1}^{-2} (\EXP{\|v_{k} - x\|^2} + \EXP{\|h_{k} - x\|^2})$ to the right-hand side to obtain for any $\zeta>0$ and $0< w_4 < 1$,
    \begin{align}
        2 \a_k^{-1} &\la F(x), \Pi_S[x_k] - x \ra \leq \a_{k-1}^{-2} \|v_{k} - x\|^2 - \a_k^{-2} \|v_{k+1} - x \|^2 + \a_{k-1}^{-2} \|h_k - x\|^2  \nonumber\\
    & - \a_k^{-2} \|h_{k+1} - x\|^2 + (\a_k^{-2} - \a_{k-1}^{-2}) \left( \|v_{k} - x\|^2 + \|h_k - x\|^2 \right) \nonumber\\
    & - \frac{\b(2-\beta)}{M_g^2 \a_k^{2}} \,\sum_{i=1}^{N_{k}} (g_{\o_{k}^{i}}^+(z_{k}^{i-1}))^2 + \frac{\zeta}{\a_k^{2}} \dist^2(x_k, S) + 2 \a_k^{-1} \la b_{k+1}^2, u_{k+1} - h_k \ra  \nonumber\\
    &+ \left[ \left( \frac{1}{\zeta} + \frac{1}{w_4} \right) \|F(x)\|^2 + 6M^2 + 7 \|b_{k+1}^2\|^2 + 6 \|b_{k+1}^1\|^2 \right]  . \label{inv_avg_eq1}
    \end{align}
    With $\a_k \leq \a_{k-1}$, the fifth quantity on the right-hand side of relation~\eqref{inv_avg_eq1} can be upper estimated using Assumption~\ref{asum_bounded_Y} as $(\a_k^{-2} - \a_{k-1}^{-2}) \left( \|v_{k} - x\|^2 + \|h_k - x\|^2 \right) \leq 2D (\a_k^{-2} - \a_{k-1}^{-2})$.
    Substituting the preceding relation back in relation~\eqref{inv_avg_eq1} and summing its both sides from $k = 1$ to $T$, we obtain
\begin{align}
    2 \sum_{k=1}^T \a_k^{-1} &\la F(x), \Pi_S[x_k] - x \ra \leq \a_0^{-2} \|v_{1} - x\|^2 - \a_T^{-2} \|v_{T+1} - x \|^2 + \a_{0}^{-2} \|h_1 - x\|^2  \nonumber\\
    & - \a_T^{-2} \|h_{T+1} - x\|^2 + 2D (\a_T^{-2} - \a_{0}^{-2}) - \frac{\b(2-\beta)}{M_g^2} \sum_{k=1}^T \a_k^{-2} \sum_{i=1}^{N_{k}} (g_{\o_{k}^{i}}^+(z_{k}^{i-1}))^2 \nonumber\\
    & + {\zeta} \sum_{k=1}^T \a_k^{-2} \dist^2(x_k, S) + 2 \sum_{k=1}^T \a_k^{-1} \la b_{k+1}^2, u_{k+1} - h_k \ra  \nonumber\\
    &+ \left( \frac{1}{\zeta} + \frac{1}{w_4} \right) \|F(x)\|^2 T + 6M^2 T + 7 \sum_{k=1}^T \|b_{k+1}^2\|^2 + 6 \sum_{k=1}^T \|b_{k+1}^1\|^2 . \label{inv_avg_eq2}
\end{align}
With $\widetilde x_T = \frac{\sum_{k=1}^T \a_k^{-1} \Pi_S[x_k]}{\sum_{k=1}^T \a_k^{-1}}$, the left-hand side of relation~\eqref{inv_avg_eq2} can be written as $\sum_{k=1}^T \a_k^{-1} \la F(x), \Pi_S[x_k] - x \ra = \left(  \sum_{k=1}^T \a_k^{-1} \right) \la F(x), \widetilde x_T - x \ra$. Substituting this into relation~\eqref{inv_avg_eq2}, taking the maximum over $x \in S$, discarding the second and fourth non-positive terms, and using Assumption~\ref{asum_bounded_Y} to upper-bound the first and third terms on the right-hand side yields the following relation
\begin{align}
    &2 \left(  \sum_{k=1}^T \a_k^{-1} \right) G(\widetilde x_T) \leq 2D \a_T^{-2} - \frac{\b(2-\beta)}{M_g^2} \sum_{k=1}^T \a_k^{-2} \sum_{i=1}^{N_{k}} (g_{\o_{k}^{i}}^+(z_{k}^{i-1}))^2 \nonumber\\
    & + {\zeta} \sum_{k=1}^T \a_k^{-2} \dist^2(x_k, S) + 2 \sum_{k=1}^T \a_k^{-1} \la b_{k+1}^2, u_{k+1} - h_k \ra  \nonumber\\
    &+ \left( \frac{1}{\zeta} + \frac{1}{w_4} \right) \left(\max_{x \in S}\|F(x)\|^2 \right) T + 6M^2 T + 7 \sum_{k=1}^T \|b_{k+1}^2\|^2 + 6 \sum_{k=1}^T \|b_{k+1}^1\|^2 , \nonumber
\end{align}
where the dual gap function $G(\cdot)$ is defined in relation~\eqref{dg_gap}. Next, we take expectation on both sides of the preceding relation, apply Assumption~\ref{asum_stoc_error}, and follow a similar analysis as in Lemma~\ref{lem_Kor_normal_avg}. The expected value of the fourth term on the right-hand side vanishes. Additionally, from~\eqref{bounded_map_eq}
we have $\max_{x \in S} \|F(x)\|^2 \leq B$. We also group the expected values of the second and third terms on the right-hand side, apply the second expression from Theorem~\ref{thm_inf_updates}, and obtain the following relation
\begin{align}
    \EXP{G(\widetilde x_T)} \leq 
    &\frac{D \a_T^{-2}}{\sum_{k=1}^T \a_k^{-1}} - \frac{\EXP{\sum_{k=1}^T \a_k^{-2} ((1-q)^{-N_k} - (1+\zeta)) \dist^2(x_k,S)}}{2 \sum_{k=1}^T \a_k^{-1}} \nonumber\\
    & + \frac{\left[ \left( \frac{1}{\zeta} + \frac{1}{w_4} \right) B^2 + 6M^2 + 13 \sigma^2 \right] T}{2 \sum_{k=1}^T \a_k^{-1}} .
\end{align}
We choose $\zeta = \frac{q}{1-q}$, which ensures $N_k \geq 1$ for all $k \geq 1$, and makes the second term on the right-hand side non-positive, so it can be dropped. Moreover, applying relation~\eqref{inf_rel} with $\gamma_k = \alpha_k^{-1}$ for all $1 \leq k \leq T$ and $\widehat x_T = \frac{\sum_{k=1}^T \alpha_k^{-1} x_k}{\sum_{k=1}^T \alpha_k^{-1}}$, we obtain
\begin{align*}
    \EXP{\widehat G(\widehat x_T)} \leq \frac{D \a_T^{-2}}{\sum_{k=1}^T \a_k^{-1}} + \frac{\left[ \left( \frac{1-q}{q} + \frac{1}{w_4} \right) B^2 + 6M^2 + 13 \sigma^2 \right] T}{2 \sum_{k=1}^T \a_k^{-1}} + \EXP{\text{V}(\widehat x_T)} .
\end{align*}
Choosing $\a_k = \min \left( \frac{\bar \alpha}{\sqrt{k+1}}, \frac{\sqrt{1 - w_4}}{\sqrt{2}L} \right)$ and applying Lemma~\ref{lem_seq} to upper bound the first two quantities yields the desired relation.
\end{proof}
Theorem~\ref{thm_Kor2} establishes that $\EXP{\widehat G(\widehat x_T)}$ converges at rate $O(1/\sqrt{T})$ along with an additional infeasibility term. The next remark, analogous to Remark~\ref{rem_sample_root_growth}, characterizes the decay rate of this infeasibility term.
\begin{remark}\label{rem_infeas_conv_inv_avg}
    For any $N_k \geq 1$ for $1 \leq k \leq T$, Lemma~\ref{lem_inf_gap_bd}(i) and Remark~\ref{rem_sample_root_growth} show $\EXP{\text{V}(\widehat x_T)} \leq B \sqrt{D} \max\limits_{1\leq k \leq T} (1-q)^{\frac{N_k}{2}}$, i.e., a geometric decay with respect to $\min_{1 \leq k \leq T} N_k$. When $N_k = \lceil k^{\frac{1}{r}} \rceil$ with $r \geq 1$, and we choose $\gamma_k = \a_k^{-1}$ in Lemma~\ref{lem_inf_gap_bd}(ii) with $\a_k = \min \left( \frac{\bar \alpha}{\sqrt{k+1}}, \frac{\sqrt{1 - w_4}}{\sqrt{2}L} \right)$ for all $1 \leq k \leq T$, then $\max_{1 \leq k \leq T} \alpha_k^{-1} = \max \left( \frac{\sqrt{T+1}}{\bar \alpha} , \frac{\sqrt{2} L}{\sqrt{1-w_4}} \right)$ and $\sum_{k=1}^T \alpha_k^{-1} \geq \max \left( \frac{1}{\bar \alpha} \left[ \left( \frac{3}{2} \right)^{\frac{1}{2}} - \frac{2}{3} \right] T^{\frac{3}{2}}, \frac{\sqrt{2} L}{\sqrt{1-w_4}} T \right)$. These relations when combined with Lemma~\ref{lem_inf_gap_bd}(ii) yield
    \begin{align*}
        \EXP{\mathrm{V}(\widehat x_T)} \leq B \sqrt{D} \left(\frac{\max \left( \frac{\sqrt{T+1}}{\bar \alpha} , \frac{\sqrt{2} L}{\sqrt{1-w_4}} \right) \left( 2^r \Gamma(r+1) \right)}{\max \left( \frac{1}{\bar \alpha} \left[ \left( \frac{3}{2} \right)^{\frac{1}{2}} - \frac{2}{3} \right] T^{\frac{3}{2}}, \frac{\sqrt{2} L}{\sqrt{1-w_4}} T \right) \left(\ln \left(\frac{1}{1-q} \right)\right)^r} \right) .
    \end{align*}
    Hence, for large $T$, $\EXP{\mathrm{V}(\widehat x_T)}$ decays at rate $O(1/T)$. Consequently, by Theorem~\ref{thm_Kor2}, $\EXP{\widehat G(\widehat x_T)}$ converges at rate $O(1/\sqrt{T})$ always.
\end{remark}
\begin{remark}\label{rem_Kor_parameter_free}
    The step-size choices in Theorems~\ref{thm_Kor1} and \ref{thm_Kor2} are upper-bounded by $\frac{\sqrt{1 - w_4}}{\sqrt{2}L}$ with $0<w_4<1$. Under Assumption~\ref{asum_bounded_Y}, this cap can be removed with minor changes to the analysis while preserving the same rates. Near the end of Lemma~\ref{lem_Kor1}, the term $\frac{2 L^2 \a_k^2}{w_2} \|x_k-u_{k+1}\|^2$ can be bounded via $\|x_k-u_{k+1}\|^2 \leq D$, giving $\frac{2 L^2 D \a_k^2}{w_2}$. Using this bound in the proof of Lemma~\ref{lem_Kor1} removes the need for a fixed upper bound on the step size; the inequality merely gains additional constants multiplying $\a_k^2$. With a decreasing step size schedule $\a_k = \frac{\bar \alpha}{\sqrt{k+1}}$, all previous results and convergence rates still hold. Moreover, the step size rule becomes problem parameter free.
\end{remark}


\section{Modified Stochastic Popov Method} \label{sec:Popov_method}

The Korpelevich method (Algorithm~\ref{algo_Kor_method}) requires two stochastic mapping evaluations per iteration, which can be prohibitive when each evaluation is expensive. To mitigate this, we use the Popov method~\cite{popov1980modification} for \eqref{VI_problem}, which needs only one new evaluation per iteration. We modify the standard Popov method by incorporating feasibility updates (Algorithm~\ref{algo_proj_steps}), yielding Algorithm~\ref{algo_Popov_method} with two update steps. The first step generates an auxiliary iterate $u_k$ using the old stochastic mapping value $\widehat F(u_{k-1}, \xi_{k-1})$, for $\xi_{k-1} \in \Xi$, and the second step computes a new stochastic mapping $\widehat F(u_k, \xi_k)$ at the auxiliary iterate $u_k$, for any $\xi_k \in \Xi$ and all $k \geq 1$. The method is initialized at some random point $x_0 \in Y$, and we choose $u_0 = x_0$. 

%
%
%
\begin{algorithm}
		\caption{Modified Stochastic Popov Method}
		\label{algo_Popov_method}
		\begin{algorithmic}[1]
			\REQUIRE{ Initial iterate $x_{0}$, step size $\a_k$}
                \FOR{$k=1,\ldots$}
                \STATE \textbf{Use} old mapping and \textbf{update:} $u_{k} = \Pi_{Y} \left[x_{k-1} - \a_{k-1} \widehat F(u_{k-1}, \xi_{k-1}) \right]$
                \STATE \textbf{Sample} $\xi_k \in \Xi$ and \textbf{update:} $v_{k} = \Pi_{Y} \left[x_{k-1} - \a_{k-1} \widehat F(u_{k}, \xi_k) \right]$
                \STATE \textbf{Call Randomized Feasibility Algorithm \ref{algo_proj_steps}:} Pass $v_{k}$ and $N_{k}$ and get $x_{k}$
                \ENDFOR
		\end{algorithmic}
\end{algorithm}	

Next, we analyze the convergence of Algorithm~\ref{algo_Popov_method}. Compared with the modified Korpelevich method (Algorithm~\ref{algo_Kor_method}), the sigma-algebra $\cF_k$ (cf. relation~\eqref{sigma_algebra}) requires refinement in this setting. Specifically with $\mathcal{F}_0 = \{x_{0}\}$, we define the sigma algebra as
\begin{align*} 
      \mathcal{F}_{k} = \mathcal{F}_0 \cup \{\omega_{t}^{i} \mid 1 \leq i \leq N_{t}, 1 \leq t \leq k\} \cup \{\xi_t \mid 1 \leq t \leq k+1 \} \quad \text{for all $k \geq 1$} . 
\end{align*}
With the defined sigma-algebra, the proofs of the feasibility update algorithm in Theorem~\ref{thm_inf_updates} and Lemma~\ref{lem_inf_geom} remain valid for Algorithm~\ref{algo_Popov_method}. Analogous to the modified Korpelevich method, we define the stochastic error for the modified Popov method for any stochastic sample $\xi_k \in \Xi$ at iteration $k \geq 1$ as
\begin{align}
    \bar b_k = F(u_k) - \widehat F(u_k, \xi_k) . \label{error_Popov}
\end{align}

We now state a lemma that underpins the convergence rate analysis of this method.
\begin{lemma}\label{lem_Popov1}
    Let Assumptions~\ref{asum_closed_set} and \ref{asum_monotone} hold and $\a_k > 0$. Then, the following relation holds for the iterates of Algorithm~\ref{algo_Popov_method} for any $x \in Y$,
    \begin{align*}
        &2 \a_k \la F(x), \Pi_S[x_k] - x \ra \leq \| x_{k} - x \|^2 - \| v_{k+1} - x \|^2 - (1-w_1) \| v_{k+1} - u_{k+1} \|^2 \nonumber\\
        & - (1-w_4) \| x_{k} - u_{k+1} \|^2 + 2 \a_k \la \bar b_{k+1}, u_{k+1} - x \ra  + \left( \frac{1}{\zeta} + \frac{1}{w_4} \right) \|F(x)\|^2 \a_k^2 \nonumber\\
        & + \left[ \frac{2 \|\bar b_{k+1}\|^2}{w_1} + \frac{2 \|\bar b_{k}\|^2}{w_1} \right] \a_k^2 + 2\a_k \|F(u_{k}) - F(u_{k+1})\| \|v_{k+1} - u_{k+1}\| + \zeta \dist^2(x_k, S) ,
    \end{align*}
    where $w_1 , w_4, \zeta > 0$ are some constants and $\bar b_k$ is defined in relation~\eqref{error_Popov}.
\end{lemma}
\begin{proof}
    Using the non-expansiveness property of the projection on the definition of $v_{k+1}$ from Algorithm~\ref{algo_Popov_method}, we obtain the following relation for any $x \in Y$
    \begin{align}
        \| v_{k+1} - x \|^2 \leq \| x_{k} - x \|^2 - \| v_{k+1} - x_{k} \|^2 + 2 \a_k \la \widehat F(u_{k+1}, \xi_{k+1}), x - v_{k+1} \ra . \label{Pop_bound2}
    \end{align}
    The second term on the right-hand side of relation~\eqref{Pop_bound2} can be written as
    \begin{align}
        \| v_{k+1} - x_{k} \|^2 = \| v_{k+1} - u_{k+1} \|^2 + \| x_{k} - u_{k+1} \|^2 - 2 \la x_{k} - u_{k+1}, v_{k+1} - u_{k+1} \ra . \label{pop_secterm}
    \end{align}
    The third term on the right-hand side of relation~\eqref{Pop_bound2} can be written as
    \begin{align}
        &\la \widehat F(u_{k+1}, \xi_{k+1}), x - v_{k+1} \ra = \la F(u_{k+1}), x - u_{k+1} \ra + \la \bar b_{k+1}, u_{k+1} - x \ra \nonumber\\
        &+ \la \widehat F(u_{k}, \xi_{k}), u_{k+1} - v_{k+1} \ra + \la \widehat F(u_{k+1}, \xi_{k+1}) - \widehat F(u_{k}, \xi_{k}), u_{k+1} - v_{k+1} \ra, \nonumber
    \end{align}
    where the quantity $\bar b_{k+1}$ is defined in relation~\eqref{error_Popov}.
    By substituting the preceding relation and relation~\eqref{pop_secterm} back in relation~\eqref{Pop_bound2}, we obtain
    \begin{align}
        \| v_{k+1} - x \|^2 &\leq \| x_{k} - x \|^2 - \| v_{k+1} - u_{k+1} \|^2 - \| x_{k} - u_{k+1} \|^2 + 2 \a_k \la \bar b_{k+1}, u_{k+1} - x \ra \nonumber\\
        & + 2 \la x_{k} - \a_k F(u_{k},\xi_k) - u_{k+1}, v_{k+1} - u_{k+1} \ra + 2 \a_k \la F(u_{k+1}), x - u_{k+1} \ra \nonumber\\
        & + 2\a_k \la \widehat F(u_{k}, \xi_{k}) - \widehat F(u_{k+1}, \xi_{k+1}), v_{k+1} - u_{k+1} \ra . \label{Pop_no1} 
    \end{align}
    Using the definition of $u_{k+1}$ from Algorithm~\ref{algo_Popov_method}, the fact that $v_{k+1} \in Y$, and the projection properties, the fifth term on the right-hand side can be dropped.
    The sixth term on the right-hand side of relation~\eqref{Pop_no1} can be upper estimated using Assumption~\ref{asum_monotone} and Young's inequality similar to relation~\eqref{split_for_dg} in Lemma~\ref{lem_Kor2} to obtain
    \begin{align}
        2 \a_k \la F(u_{k+1}), x - u_{k+1} \ra \leq &2 \a_k \la F(x), x - \Pi_S[x_k] \ra + \left[ \frac{1}{\zeta} + \frac{1}{w_4} \right] \|F(x)\|^2 \a_k^2 \nonumber\\
        &+ \zeta \dist^2(x_k, S) + w_4 \|u_{k+1}-x_k\|^2 , \label{split_for_dg_popov}
    \end{align}
    where $\zeta$ and $w_4$ are positive constants. The last term on the right-hand side of relation~\eqref{Pop_no1}
    can be written as
    \begin{align}
        2\a_k \la \widehat F(u_{k}, \xi_{k}) &- \widehat F(u_{k+1}, \xi_{k+1}), v_{k+1} - u_{k+1} \ra = 2\a_k \la F(u_{k}) - F(u_{k+1}), v_{k+1} - u_{k+1} \ra \nonumber\\
        &+ 2\a_k \la \bar b_{k+1}, v_{k+1} - u_{k+1} \ra - 2\a_k \la \bar b_{k}, v_{k+1} - u_{k+1} \ra . \nonumber
    \end{align}
    Applying Cauchy-Schwarz inequality, the growth condition on the mapping $F$ (Assumption~\ref{asum_map_growth}), and Young's inequality, we obtain for any $w_1 > 0$,
    \begin{align}
        &2\a_k \la \widehat F(u_{k}, \xi_{k}) - \widehat F(u_{k+1}, \xi_{k+1}), v_{k+1} - u_{k+1} \ra \nonumber\\
        &\leq 2\a_k (\|\bar b_{k+1}\| + \|\bar b_{k}\| + \|F(u_{k}) - F(u_{k+1})\|) \|v_{k+1} - u_{k+1}\| \nonumber\\
        & \leq \frac{2 \a_k^2}{w_1} (\|\bar b_{k+1}\|^2 + \|\bar b_{k}\|^2) + w_1 \|v_{k+1} \!-\! u_{k+1}\|^2 + 2\a_k \|F(u_{k}) \!-\! F(u_{k+1})\| \|v_{k+1} \!-\! u_{k+1}\| . \nonumber
    \end{align}
    Applying the preceding relation and relation~\eqref{split_for_dg_popov} in relation~\eqref{Pop_no1} and dropping the fifth term on the right-hand side yields the desired relation of the lemma.
\end{proof}
The relation in Lemma~\ref{lem_Popov1} contains the term $2\a_k \|F(u_{k}) - F(u_{k+1})\| \|v_{k+1} - u_{k+1}\|$, which must be controlled. The conventional approach, applying Assumption~\ref{asum_map_growth} with Young’s inequality, introduces the term $\|u_k - u_{k+1} \|^2$ which cannot be combined with the third and fourth terms on the right-hand side of Lemma~\ref{lem_Popov1} to control the error. Instead, applying Young’s inequality directly with a constant $w_2>0$ yields
\begin{align*}
    2\a_k \|F(u_{k}) \!-\! F(u_{k+1})\| \|v_{k+1} - u_{k+1}\| \leq \frac{ \a_k^2}{w_2} \|F(u_{k}) \!-\! F(u_{k+1})\|^2 + w_2 \|v_{k+1} - u_{k+1}\|^2.
\end{align*}
Using this bound in Lemma~\ref{lem_Popov1} and setting the constants $w_1 = w_2 = \frac{1}{2}$ and $0 < w_4 < 1$, a line search can be performed to select the step size $\a_k$ so that $2 \a_k^2 \|F(u_{k}) - F(u_{k+1})\|^2 \leq (1-w_4) \|x_k - u_{k+1} \|^2$. While this removes the need for Assumption~\ref{asum_map_growth}, it adds a subproblem and computational overhead. Moreover, Lemma~\ref{lem_Popov1} still contains a term of $\a_k^2$, requiring a diminishing step size for convergence. Since line search tightens only the upper bound without improving the theoretical rate, we instead rely on Assumption~\ref{asum_bounded_Y} and, hence, the boundedness of the mapping $F$ (cf.\ relation~\eqref{bounded_map_eq}) for the convergence analysis. We next present a lemma supporting this approach.
\begin{lemma}\label{lem_Popov2}
    Under Assumptions~\ref{asum_closed_set}, \ref{asum_bounded_Y}, \ref{asum_map_growth}, and \ref{asum_monotone}, and step size $\a_k > 0$, the following relation is satisfied by the iterates of Algorithm~\ref{algo_Popov_method} for any $x \in S$,
    \begin{align*}
        2 \a_k &\la F(x), \Pi_S[x_k] - x \ra \leq \| v_{k} - x \|^2 - \| v_{k+1} - x \|^2 + \|h_k - x\|^2 - \|h_{k+1} - x\|^2 \\
        & - \frac{\beta(2-\beta)}{M_g^2} \sum_{i=1}^{N_k} (g_{\o_k^i}^+(z_k^{i-1}))^2 + \zeta \dist^2(x_k,S) + 2 \a_k \la \bar b_{k+1}, u_{k+1} - h_k \ra \\
        & + \left[ \left( \frac{1}{\zeta} + 3 \right) B^2 + 4\|\bar b_{k+1}\|^2 + 4\|\bar b_k\|^2 \right] \a_k^2 ,
    \end{align*}
    where $\bar b_k$ is defined in relation~\eqref{error_Popov} and $h_{k+1} = \Pi_S[h_k - \a_k \bar b_{k+1}]$ with $h_0 \in S$.
\end{lemma}
\begin{proof}
    We start from the relation of Lemma~\ref{lem_Popov1} and apply Assumptions~\ref{asum_map_growth} and~\ref{asum_bounded_Y} on the quantity $\|F(u_k) - F(u_{k+1})\|$ to obtain $\|F(u_k) - F(u_{k+1})\| \leq L \|u_k - u_{k+1}\| + M \leq L \sqrt{D} + M \leq B$,
    where the bound $B$ is defined in relation~\eqref{bounded_map_eq}. Applying relation~\eqref{bounded_map_eq} and the preceding bound to upper estimate the sixth and the eighth terms on the right-hand side of the relation of Lemma~\ref{lem_Popov1}, we obtain
    \begin{align*}
        &2 \a_k \la F(x), \Pi_S[x_k] - x \ra \leq \| x_{k} - x \|^2 - \| v_{k+1} - x \|^2 - (1-w_1) \| v_{k+1} - u_{k+1} \|^2 \nonumber\\
        & - (1-w_4) \| x_{k} - u_{k+1} \|^2 + 2 \a_k \la \bar b_{k+1}, u_{k+1} - x \ra  + \left( \frac{1}{\zeta} + \frac{1}{w_4} \right) B^2 \a_k^2 \nonumber\\
        & + \left[ \frac{2 \|\bar b_{k+1}\|^2}{w_1} + \frac{2 \|\bar b_{k}\|^2}{w_1} \right] \a_k^2 + 2\a_k B \|v_{k+1} - u_{k+1}\| + \zeta \dist^2(x_k, S) .
    \end{align*}
    The eigth term on the right-hand side can be upper estimated using Young's inequality, $2\a_k B \|v_{k+1} - u_{k+1}\| \leq \frac{B^2 \a_k^2}{w_2} + w_2 \|v_{k+1} - u_{k+1}\|^2$,
    which when substituted back in the preceding relation along with the constants $w_1 = w_2 = \frac{1}{2}$ and $w_4 = 1$ yields
    \begin{align*}
        &2 \a_k \la F(x), \Pi_S[x_k] - x \ra \leq \| x_{k} - x \|^2 - \| v_{k+1} - x \|^2 + 2 \a_k \la \bar b_{k+1}, u_{k+1} - x \ra  \nonumber\\
        & + \left( \frac{1}{\zeta} + 3 \right) B^2 \a_k^2 + 4 \left( \|\bar b_{k+1}\|^2 + \|\bar b_{k}\|^2 \right) \a_k^2 + \zeta \dist^2(x_k, S) .
    \end{align*}
    The rest of the analysis follows along the same lines as the proof of Lemma~\ref{lem_Kor2}.
\end{proof}
The relation in Lemma~\ref{lem_Popov2} is analogous to Lemma~\ref{lem_Kor2} for the modified Korpelevich method. Consequently, the convergence rate analysis for the modified Popov method proceeds along similar lines. Without a detailed proof, we therefore state the following theorem, which provides an upper bound on the function $\widehat G(\cdot)$ (cf. relation~\eqref{mod_dual_gp}) for Algorithm~\ref{algo_Popov_method} under a diminishing step size $\a_k = \frac{\bar \alpha}{\sqrt{t+1}}$. 
\begin{theorem}\label{thm_Popov}
    Under Assumptions~\ref{asum_closed_set}, \ref{asum_bounded_Y}, \ref{asum_map_growth},  \ref{asum_monotone}, \ref{asum-regularmod}, and \ref{asum_stoc_error}, the following hold for Algorithm~\ref{algo_Popov_method} with feasibility step sizes $0 < \beta < 2$ and constant $\bar \alpha>0$.\\
    i) For $\alpha_k = \frac{\bar \alpha}{\sqrt{k+1}}$, $k \geq 1$, the averaged iterate $\widehat x_T = \frac{\sum_{k=1}^T \a_k x_k}{\sum_{k=1}^T \a_k}$ satisfies for any $T \geq 1$,
    \begin{align*}
        \EXP{\widehat G(\widehat x_T)} &\leq \frac{\frac{\sqrt{2} D}{\bar \alpha} + \frac{\bar \alpha}{\sqrt{2}} \left[ \left( \frac{1-q}{q} + 3 \right) B^2 + 8 \sigma^2 \right] \ln(T+1)}{\sqrt{T}} + \EXP{\mathrm{V}(\widehat x_T)} ,
    \end{align*}
    where $\Gamma(r) = \int_{0}^\infty p^{r-1} \exp(-p) dp$ is the gamma function. The constants $q$, $B$ and $D$ are defined in relations~\eqref{quantity-q}, \eqref{bounded_map_eq}, and Assumption~\ref{asum_bounded_Y}, respectively, while $\mathrm{V}(\widehat x_T)$ is an infeasibility violation term defined in relation~\eqref{inf_rel}.

    When $\alpha_k = \frac{\bar \alpha}{\sqrt{T}}$, $k \geq 1$, and $\widehat x_T = \frac{1}{T} \sum_{k=1}^T x_k$, the following holds for any $T \geq 1$,
    \begin{align*}
        \EXP{\widehat G(\widehat x_T)} &\leq \frac{\frac{D}{\bar \alpha} + \bar \alpha \left[ \left( \frac{1-q}{q} + 3 \right) B^2 + 8 \sigma^2 \right]}{\sqrt{T}} + \EXP{\mathrm{V}(\widehat x_T)} .
    \end{align*}
    ii) For $\alpha_k = \frac{\bar \alpha}{\sqrt{k+1}}$, $k \geq 1$, the average  
    $\widehat x_T = \frac{\sum_{k=1}^T \a_k^{-1} x_k}{\sum_{k=1}^T \a_k^{-1}}$ satisfies for any $T \geq 2$,
    \begin{align*}
        \EXP{\widehat G(\widehat x_T)} &\leq \frac{\frac{3 D}{\bar \alpha} + \bar \alpha \left[ \left( \frac{1-q}{q} + 3 \right) B^2 + 8 \sigma^2 \right]}{2 \left[ \left( \frac{3}{2} \right)^{\frac{1}{2}} - \frac{2}{3} \right] \sqrt{T}} + \EXP{\mathrm{V}(\widehat x_T)} .
    \end{align*}
\end{theorem}
\begin{proof}
    \textit{(i):} With $\widehat x_T = \frac{\sum_{k=1}^T \a_k x_k}{\sum_{k=1}^T \a_k}$, we start from Lemma~\ref{lem_Popov2} and follow the same analysis as done for Lemma~\ref{lem_Kor_normal_avg}, yielding for $T \geq 1$,
\begin{align*}
    \EXP{\widehat G(\widehat x_T)} &\leq \frac{D}{\sum_{k=1}^T \a_k} + \frac{\left[ \left( \frac{1-q}{q} + 3 \right) B^2 + 8 \sigma^2 \right] \sum_{k=1}^T \a_k^2}{2 \sum_{k=1}^T \a_k} + \EXP{\mathrm{V}(\widehat x_T)} . 
\end{align*}
For $\alpha_k = \frac{\bar \alpha}{\sqrt{k+1}}$, the first two quantities on the right-hand side can be upper estimated using Lemma~\ref{lem_seq}, leading to the first relation of part (i). The second relation of part (i) can be obtained similarly.

\noindent \textit{(ii):} Here, $\widehat x_T = \frac{\sum_{k=1}^T \a_k^{-1} x_k}{\sum_{k=1}^T \a_k^{-1}}$. Starting from Lemma~\ref{lem_Popov2}, and following the steps of Theorem~\ref{thm_Kor2}, we obtain
\begin{align*}
    \EXP{\widehat G(\widehat x_T)} &\leq \frac{D \a_T^{-2}}{\sum_{k=1}^T \a_k^{-1}} + \frac{\left[ \left( \frac{1-q}{q} + 3 \right) B^2 + 8 \sigma^2 \right] T}{2 \sum_{k=1}^T \a_k^{-1}} + \EXP{\text{V}(\widehat x_T)} .
\end{align*}
With $\alpha_k = \frac{\bar \alpha}{\sqrt{k+1}}$ for all $k \geq 1$, the last relation of Lemma~\ref{lem_seq} yields the result.
\end{proof}
Theorem~\ref{thm_Popov} stated for Algorithm~\ref{algo_Popov_method} is analogous to the Theorems~\ref{thm_Kor1} and \ref{thm_Kor2} stated for Algorithm~\ref{algo_Kor_method}. Following Lemma~\ref{lem_inf_gap_bd}(i), $\EXP{\text{V}(\widehat x_k)} \leq B \sqrt{D} \left(\max_{1 \leq k \leq T} (1-q)^{\frac{N_k}{2}}  \right)$ for any $\{N_k\}_{k=1}^T$. Moreover, if we choose the sample size $N_k = \lceil k^{\frac{1}{r}} \rceil$, $r \geq 1$, for $1 \leq k \leq T$, then by similar analysis as in Remarks~\ref{rem_sample_root_growth} and \ref{rem_infeas_conv_inv_avg}, we obtain
\begin{align*}
    \EXP{\text{V}(\widehat x_T)} \leq \begin{cases}
            \frac{B \sqrt{D}}{\sqrt{T}} \frac{2^r \Gamma(r+1)}{\left(\ln \left(\frac{1}{1-q} \right)\right)^r} &\text{if $\alpha_k = \frac{\bar \alpha}{\sqrt{k+1}}$, $\widehat x_T = \frac{\sum_{k=1}^T \a_k x_k}{\sum_{k=1}^T \a_k}$, $T \geq 1$,}\\
            \frac{B \sqrt{D}}{T} \frac{2^r \Gamma(r+1)}{\left(\ln \left(\frac{1}{1-q} \right)\right)^r} &\text{if $\alpha_k = \frac{\bar \alpha}{\sqrt{T}}$, $\widehat x_T = \frac{1}{T} \sum_{k=1}^T x_k$, $T \geq 1$,}\\
            \frac{\sqrt{3} B \sqrt{D}}{\sqrt{2} \left[ \left( \frac{3}{2} \right)^{\frac{1}{2}} - \frac{2}{3} \right] T} \frac{2^r \Gamma(r+1)}{\left(\ln \left(\frac{1}{1-q} \right)\right)^r} &\text{if $\alpha_k = \frac{\bar \alpha}{\sqrt{k+1}}$, $\widehat x_T = \frac{\sum_{k=1}^T \a_k^{-1} x_k}{\sum_{k=1}^T \a_k^{-1}}$, $T \geq 2$.}
        \end{cases}
\end{align*}
For $N_k = \lceil k^{\frac{1}{r}} \rceil$, the preceding relation and Theorem~\ref{thm_Popov} show a best case convergence rate of $O(1/\sqrt{T})$ both with a constant step size and with a diminishing step size with $\{\a_k^{-1}\}_{k=1}^T$ as averaging weights. In this case, the infeasibility gap decays as $O(1/T)$. By contrast, under the first statement of the theorem, the infeasibility decay rate is $O(1/\sqrt{T})$, and a term that has an additional $\ln(T)$ factor.

\begin{remark}
    Here we examine the convergence behavior of the infeasibility term when $N_k = \lceil \log_m(k+1) \rceil$ with base $m > 1$. Following the analysis of Lemma~\ref{lem_inf_gap_bd}(ii), we can upper estimate $\sum_{k=1}^T \left(\sqrt{1-q} \right)^{N_k} \leq \int_{0}^T \exp(- a \ln(k+1)) dk = \int_{0}^T (k+1)^{-a} dk$, where $a = \frac{\ln \left( \frac{1}{1-q} \right)}{2 \ln(m)} = \frac{1}{2} \log_m \left( \frac{1}{1-q} \right)$. The preceding integration simplifies for all $T \geq 1$,
    \begin{align*}
        \sum_{k=1}^T \left(\sqrt{1-q} \right)^{\lceil \log_m(k+1) \rceil} 
        \leq 
        \begin{cases}
            \frac{1-(T+1)^{-(a-1)}}{a-1} \leq \frac{1}{a-1} & \text{if $a > 1$, i.e., $q > 1 - m^{-2}$,}\\
            \frac{(T+1)^{1-a} - 1}{1-a} \leq \frac{(T+1)^{1-a}}{1-a} & \text{if $a < 1$, i.e., $q < 1 - m^{-2}$,} \\
            \ln(T+1) & \text{if $a = 1$, i.e., $q = 1 - m^{-2}$.}
        \end{cases}
    \end{align*}
    If $a \leq 1$ (i.e., $q \leq 1 - m^{-2}$) the infeasibility error grows. Using the analysis of Lemma~\ref{lem_inf_gap_bd}(ii) and Theorem~\ref{thm_Popov}, with either a constant step size $\alpha_k = \frac{\bar \alpha}{\sqrt{T}}$ or a diminishing step size $\alpha_k = \frac{\bar \alpha}{\sqrt{k+1}}$ and averaged iterate $\widehat x_k = \frac{\sum_{k=1}^T \a_k^{-1} x_k}{\sum_{k=1}^T \a_k^{-1}}$, the expected infeasibility term satisfies: $O(1/T)$ if $a > 1$, $O(1/T^{a})$ if $a < 1$ and $O(\ln(T)/T)$ if $a=1$. Intuitively, small steps in the first case keep the iterates near the solution set. In the second case, because $N_k$ grows with $k$, averaging with weights $\a_k^{-1}$ emphasizes later (more feasible) iterates. Similar rates also hold for the modified Korpelevich method. However, if one considers the averaged iterate $\widehat x_T$ with averaging weights $\a_k = \frac{\bar \alpha}{\sqrt{k+1}}$ for all $1 \leq k \leq T$ and $T \geq 1$ (cf. first relations of Remark~\ref{rem_sample_root_growth} and Theorem~\ref{thm_Popov}(i)), the infeasibility term may fail to converge when $1/2 <a <1$.
\end{remark}
The analysis in this work assumes a deterministic sample size $N_k$, but extends directly to stochastic $N_k$ drawn from a distribution, as studied in~\cite{chakraborty2026randomized}.



\section{Extension to Set-valued Mappings}\label{sec:multi_valued_F}

We extend the preceding algorithms and convergence results to set-valued mappings $F:Y\rightrightarrows\mathbb{R}^n$ by modifying the mapping assumptions, stochastic-oracle model, and merit function, while the random feasibility updates and their analysis remain unchanged. The set-valued Variational Inequality problem is to find $x^* \in S$ and $f^* \in F(x^*)$ such that
\begin{align}
    \la f^*, x - x^* \ra \geq 0 \qquad \text{for all } x \in S = X \cap Y . \label{prob_multi_valued_VI}
\end{align}
For every \(x,y\in Y\), \(f_x\in F(x)\), and \(f_y\in F(y)\), Assumptions~\ref{asum_map_growth} and \ref{asum_monotone} become
\begin{align*}
    \langle f_x-f_y,x-y\rangle\geq 0 \quad \text{and} \quad \|f_x-f_y\|
    \leq L\|x-y\|+M \;\; \text{for some $L, M \geq 0$.}
\end{align*}
Assumption~\ref{asum_stoc_error} is modified as
follows: there exists a measurable selection
$x\mapsto f_x\in F(x)$ such that, for every $x\in Y$,
\[
    \EXP{\widehat F(x,\xi)}=f_x,
    \qquad
    \EXP{\|\widehat F(x,\xi)-f_x\|^2}\leq\sigma^2.
\] 
Moreover, analogous to relation~\eqref{bounded_map_eq}, we see that for a fixed $x^{\mathrm{ref}}\in Y$ and
$f^{\mathrm{ref}}\in F(x^{\mathrm{ref}})$, and for every $x\in Y$ and $f_x\in F(x)$, 
\begin{align*}
    \|f_x\|
    \leq
    L\|x-x^{\mathrm{ref}}\|+M+\|f^{\mathrm{ref}}\|
    \leq
    L\sqrt{D}+M+\|f^{\mathrm{ref}}\|
    =:B_{\mathrm{mv}} .
\end{align*}
It remains to define the \textit{modified} dual gap function for the set-valued mapping case:
\[
\widehat G_{\mathrm{mv}}(y)
:=
\left|
    \sup_{x\in S}
    \sup_{f_x\in F(x)}
    \langle f_x,y-x\rangle
\right|,
\qquad \text{for all } y\in Y.
\] 
We retain the supremum because Assumption~\ref{asum_bounded_Y} (compactness of set $Y$ and hence $S$) alone does not guarantee attainment. If $\operatorname{gph}F|_S := \{(x,f)\in S\times\mathbb{R}^n:f\in F(x)\}$ is compact, for example if $F$ is upper semicontinuous on $S$ with nonempty compact values, then the supremum is attained, and both suprema can be replaced by maxima. For all $y\in S$, the dual gap function for set-valued mapping is defined as
\[
    G_{\mathrm{mv}}(y)
    :=
    \sup_{x\in S}\sup_{f_x\in F(x)}
    \langle f_x,y-x\rangle .
\]
For the averaged iterates $\widehat x_k$ and $\widetilde x_k$ defined in relation~\eqref{avg_pt}, the argument used in relations~\eqref{inf_rel} and \eqref{terminfes_error} similarly yields
\[
    \EXP{\widehat G_{\mathrm{mv}}(\widehat x_k)}
    \leq
    \EXP{G_{\mathrm{mv}}(\widetilde x_k)}
    +
    B_{\mathrm{mv}} \EXP{\|\widehat x_k-\widetilde x_k\|} .
\]
With these modifications, the convergence arguments in Sections~\ref{sec:Kor_algo} to \ref{sec:Popov_method}
carry over with only notational changes. At each oracle call, $F(x)$ is
replaced by the mean selection $f_x\in F(x)$, and the oracle error is
defined relative to this selection. The monotonicity and controlled-growth
arguments use their set-valued counterparts, while the random
feasibility results of Section~\ref{sec:rand_feas_updates} remain unchanged because they do not
involve $F$. Consequently, Theorems~\ref{thm_Kor1}, \ref{thm_Kor2},
and~\ref{thm_Popov} remain valid with $\widehat G$ and $B$ replaced by
$\widehat G_{\mathrm{mv}}$ and $B_{\mathrm{mv}}$, respectively, preserving
the same rates under the
same step size, averaging, and sample size $N_k$ choices.


\section{Simulation on a Zero-Sum Matrix Game}\label{sec:simulation}


\begin{figure*}[!h]
  \centering
  \begin{subfigure}{0.4\textwidth}
    \centering
    \includegraphics[width=\linewidth, height =\linewidth]{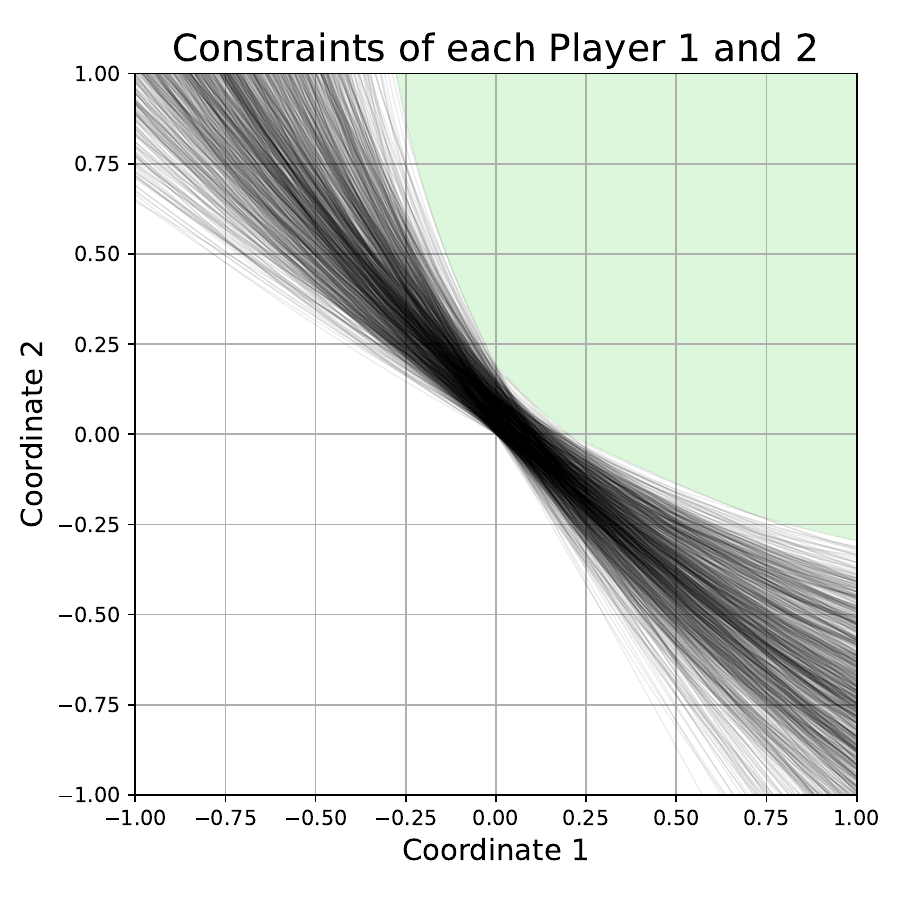}
    \caption{The feasible set of each player}\label{fig_feas_set}
  \end{subfigure}\hfill
  \begin{subfigure}{0.58\textwidth}
    \centering
    \includegraphics[width=\linewidth, height = 0.7\linewidth]{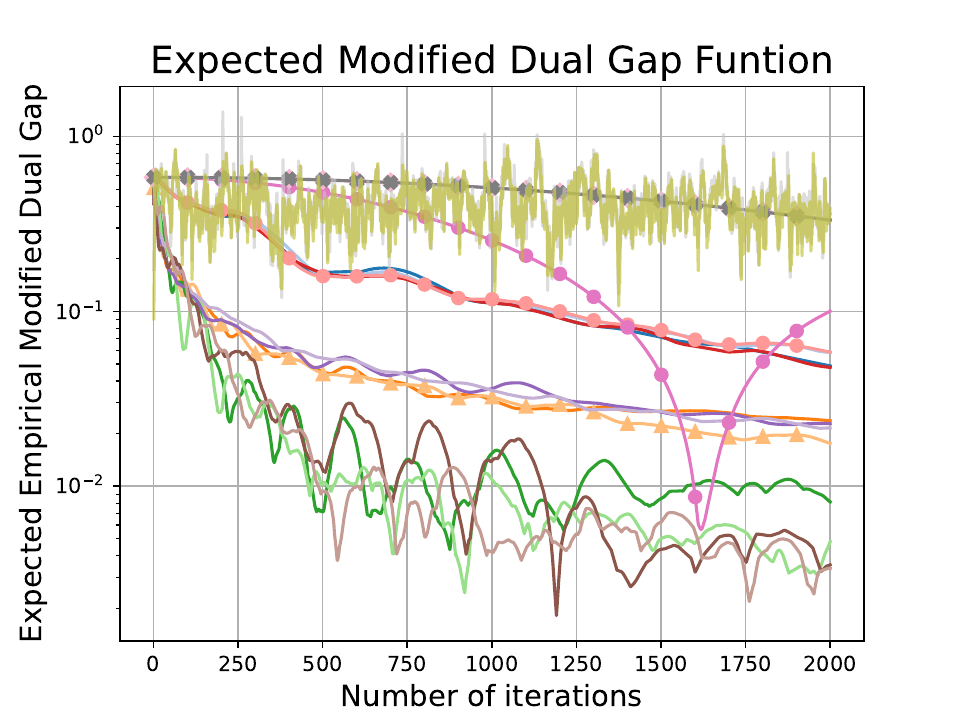}
    \caption{Modified Dual Gap Function}\label{fig_dg_gap}
    \end{subfigure}\hfill
    \begin{subfigure}{0.5\textwidth}
    \centering
    \includegraphics[width=\linewidth]{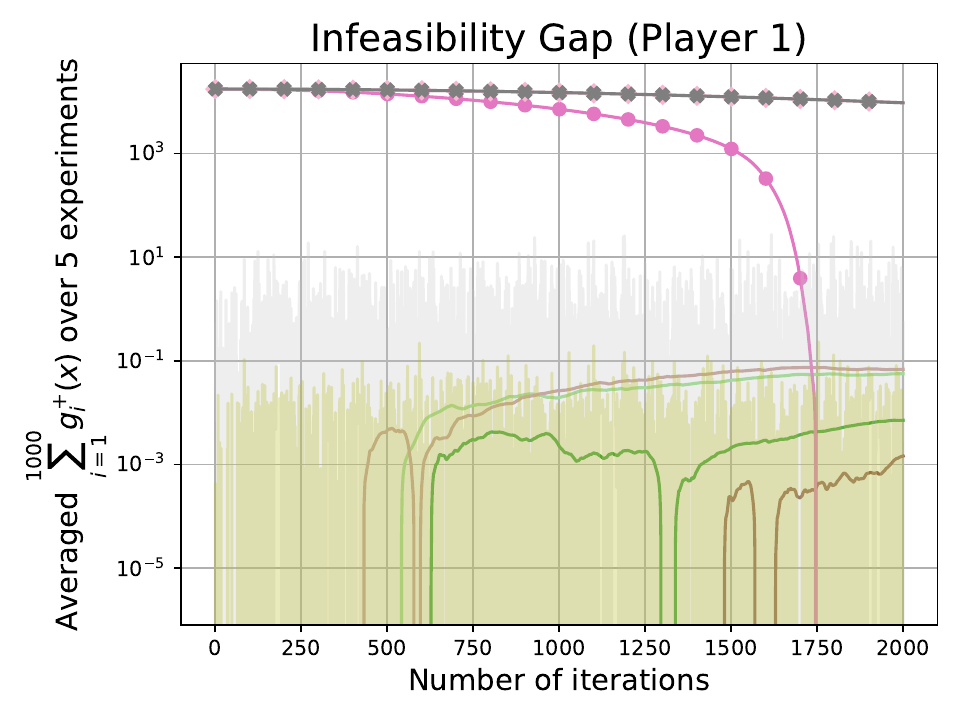}
    \caption{Infeasibility of Player 1}\label{fig_infeas_gap_ply1}
    \end{subfigure}\hfill
    \begin{subfigure}{0.5\textwidth}
    \centering
    \includegraphics[width=\linewidth]{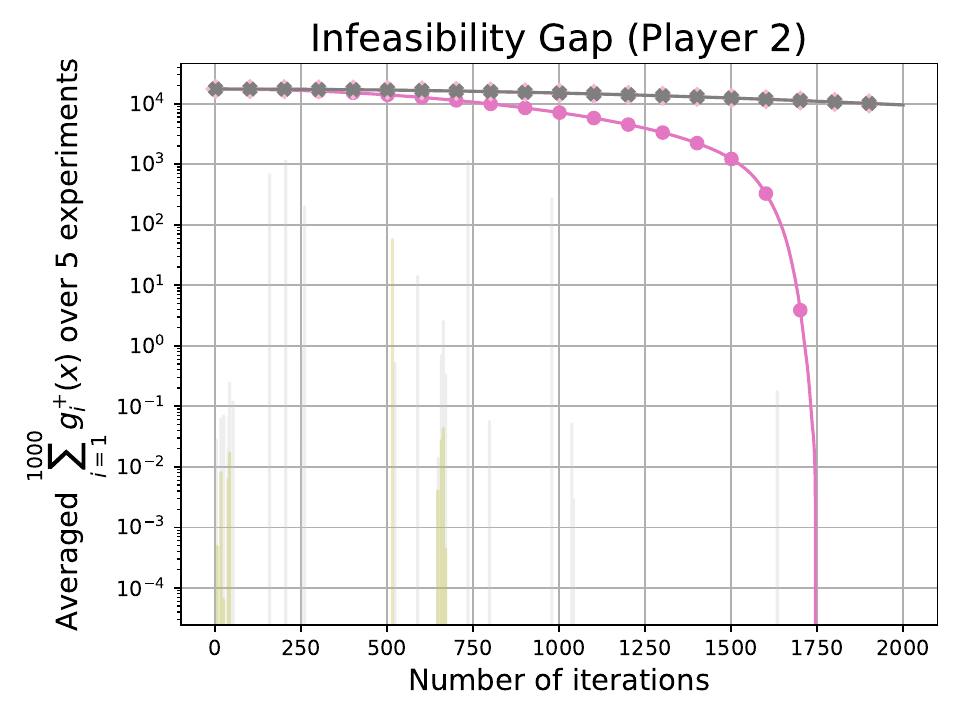}
    \caption{Infeasibility of Player 2}\label{fig_infeas_gap_ply2}
    \end{subfigure}\hfill
    \begin{subfigure}{\textwidth}
    \centering
    \includegraphics[width=\linewidth]{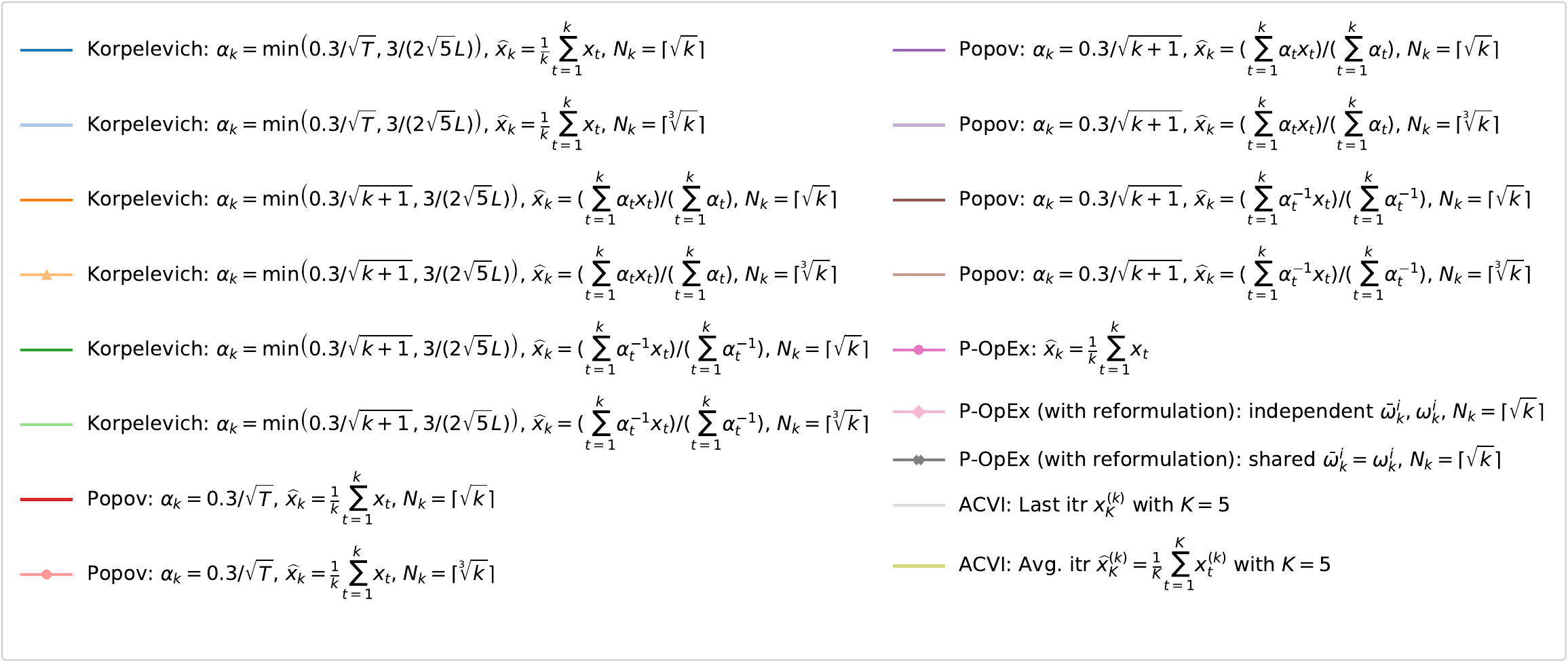}
    \caption{Legends of all the plots above}\label{fig_legend}
    \end{subfigure}\hfill

    \caption{Algorithms~\ref{algo_Kor_method} and~\ref{algo_Popov_method} are compared with P-OpEx~\cite[Algorithm~2]{boob2026first} (with and without reformulation) and ACVI~\cite[Algorithm~1]{yang2022solving} with respective averaging schemes.}
  \label{fig:overview}
\end{figure*}


We consider the constrained two-player zero-sum game formulated as
\begin{align*}
    &\min_y \max_z \EXP{y^T A z + \la \xi_1,y \ra - \la \xi_2, z \ra}, \quad \text{st. } y \in [-1 , 1]^{2} , \quad z \in [-1 , 1]^{2},\\
    &\text{and} \quad  \la y , B_i y \ra + \la c_i , y \ra - d_i \leq 0, \quad \la z , B_i z \ra + \la c_i , z \ra - d_i \leq 0, \quad i=1,\ldots, 1000 .
\end{align*}
The random variable $\xi_1, \xi_2 \sim \mathcal{N}(0,0.5^2 \mathbf{I}_2)$, where $\mathcal{N}$ denotes the normal distribution and $\mathbf{I}_2$ is the 2-dimensional identity matrix. The objective of the game reduces to $\EXP{y^T A z + \la \xi_1,y \ra - \la \xi_2, z \ra} = y^T A z$ since $\EXP{\xi_1} = \EXP{\xi_2} = 0$. Note that Assumptions~\ref{asum_closed_set} and \ref{asum_bounded_Y} are satisfied. 
Let $x = [y,z]^T$, so the SVI mapping and its expectation are
\begin{align*}
    \widehat F(x, \xi) = \begin{bmatrix}
        Az  \\
        - A^T y  
    \end{bmatrix} \!+\! \begin{bmatrix}
        \xi_1 \\
        \xi_2  
    \end{bmatrix} = \tilde A x + \xi \text{ with } \tilde A = \begin{bmatrix}
        0 & A  \\
        - A^T & 0  
    \end{bmatrix} \text{ and } \xi = \begin{bmatrix}
        \xi_1 \\
        \xi_2  
    \end{bmatrix}, \; F(x) = \tilde A x .
\end{align*}
Since $\xi \in \R^{4}$ and each element of the vector $\xi$ is sampled from $\mathcal{N}(0,0.5^2)$, Assumption~\ref{asum_stoc_error} holds with $\sigma= \sqrt{0.5^2 \times 4} = 1$. The mapping $F(\cdot)$ is monotone (Assumption~\ref{asum_monotone})  since $\langle F(x)-F(\bar x),\, x-\bar x\rangle=0$ for any $x=[y,z]$ and $\bar x=[\bar y,\bar z]$ (due to the structure of $\tilde A$). Since $\tilde A$ is skew-symmetric, it is normal (hence unitarily diagonalizable) and has complex eigenvalues. Here, the Lipschitz constant $L>0$ can be computed via the power method \cite[Algorithm 1]{gouk2021regularisation}. Assumption~\ref{asum_map_growth} holds with $M=0$ and $L=\|\tilde A\|_2=\|A\|_2$. In our experiments, we generate $A$ with eigenvalues in $[0,4]$. Specifically, we form a diagonal matrix $\Lambda$ with diagonal entries sampled i.i.d. from $U[0,4]$, draw a random matrix $M$, compute its QR factorization $M=QR$ with $Q$ orthogonal and $R$ upper triangular matrices, and set $A=Q\Lambda Q^\top$. Similarly, the matrices $B_i$ for $i\in\{1,\ldots,1000\}$ are generated with eigenvalues sampled from $U[0,2]$; the vectors $c_i$ are sampled from $U[-10,-5]$; and the scalars $d_i$ are sampled from $U[-1,0]$. The generated set is illustrated in Figure~\ref{fig_feas_set}, which shows that the constraint set has a nonempty interior; therefore, Assumption~\ref{asum-regularmod} is satisfied. 


\begin{table}[t]
\centering
\footnotesize
\caption{Run time for a single experiment of all the algorithms}
\label{tab_runtime}
\begin{tabular}{@{}ll@{}}
\toprule
\textbf{Method} & \textbf{Run time} \\
\midrule
Algorithms~\ref{algo_Kor_method} and \ref{algo_Popov_method} & $< 1$ sec \\
P-OpEx \cite[Algorithm~2]{boob2026first} & 17 secs \\
P-OpEx with reformulation (Algorithm~\ref{alg:single-expectation-popex}) & $< 1$ sec \\
ACVI \cite[Algorithm~1]{yang2022solving} & 12 mins 15 secs \\
\bottomrule
\end{tabular}
\end{table}


We compare Algorithms~\ref{algo_Kor_method} and
\ref{algo_Popov_method} with P-OpEx~\cite[Algorithm~2]{boob2026first},
which uses operator extrapolation with primal-dual updates, and
ACVI~\cite[Algorithm~1]{yang2022solving}, which applies ADMM with
log-barrier penalty functions. Standard P-OpEx considers all constraints
and associates one dual variable with each constraint, which can be
computationally expensive when the number of constraints is large.
For experimental evaluation, we therefore consider the reformulation $X = \bigcap_{a\in\mathcal{A}}\{x\mid g_a(x)\leq 0\} = \left\{
x\;\middle|\;
\mathsf{E}_{a\sim\mathcal{A}}\!\left[g_a^+(x)\right]\leq 0
\right\}$,
and run $N_k$ times per iteration P-OpEx with a single dual multiplier associated with the expectation
constraint. The detailed algorithm is provided in
Appendix~\ref{app_popex}. However, the reformulated
constraint admits no Slater point and does not satisfy MFCQ
(cf.\ Appendix~\ref{app_flaw}); consequently, no convergence theory is
currently available for the resulting method. Developing such a theory is
beyond the scope of this work, so we compare the method only empirically. All the algorithms are run on a MacBook Pro 2021 (Apple M1 pro chip, 16GB RAM). For our methods, we choose the constant and diminishing step-size schemes from Theorems~\ref{thm_Kor1}, \ref{thm_Kor2}, and \ref{thm_Popov}, with the hyperparameter set to $\bar \alpha = 0.3$ for all the algorithms and $w_4 = \frac{1}{10}$ for Algorithm~\ref{algo_Kor_method}. We consider two growth schedules for the number of inner feasibility iterations: $N_k=\lceil \sqrt{k}\rceil$ and $N_k=\lceil \sqrt[3]{k}\rceil$, with constraint indices sampled uniformly at random from $\{1,\ldots,1000\}$ with replacement. The parameters for P-OpEx are chosen according to \cite[Theorem~2]{boob2026first}: we set the square root of the diameter $D$ in Assumption~\ref{asum_bounded_Y} (corresponding to $D_X$ in \cite{boob2026first}) to $4$, compute the remaining parameters from the game objective and constraints, and let $B = 10$ in \cite{boob2026first}. Note that P-OpEx requires knowledge of all problem parameters, including the variance of the stochastic noise, which may be unavailable in practice. For ACVI, we select the hyperparameters by trial and error: $\mu_{-1}=10^{-5}$, $\delta=0.5$, $\beta=20$, and $K=5$ inner iterations.

For all algorithms, we plot the expected value of the modified dual-gap function~\eqref{mod_dual_gp} evaluated at the averaged iterates in Figure~\ref{fig_dg_gap}. The modified dual gap function is computed by empirically generating 1500 points in the box $[-1,1]^2$ for both players and retaining only the feasible points, i.e., any point $\bar y \in [-1,1]^2$ for which $\la \bar y , B_i \bar y \ra + \la c_i , \bar y \ra - d_i \leq 0$. The expectation is estimated by averaging over five independent runs. As a proxy for infeasibility for both players, we report the average of $\sum_{i=1}^{1000} g_i^{+}(x)$ over the same five different experimental runs (cf.\ Figures~\ref{fig_infeas_gap_ply1} and \ref{fig_infeas_gap_ply2}), where $x$ denotes the decision block per player of the averaged iterates produced by each method. We use this surrogate since the exact distance from an iterate to the feasible set $S$ (via projection) is not available. Figure~\ref{fig_dg_gap} shows that Algorithms~\ref{algo_Kor_method} and~\ref{algo_Popov_method} perform better in terms of the empirically expected modified dual gap function than P-OpEx (both with and without constraint reformulation) and ACVI, while also requiring much less computational time (see Table~\ref{tab_runtime}). From Figure~\ref{fig_infeas_gap_ply1}, it is evident that infeasibility reduces with a larger sample-size, i.e., $N_k = \lceil \sqrt{k} \rceil$ outperforms those with $N_k = \lceil \sqrt[3]{k} \rceil$. Moreover, Algorithms~\ref{algo_Kor_method} and \ref{algo_Popov_method} with $\alpha_k^{-1}$-weighted averaging perform best in
Figure~\ref{fig_dg_gap} under the stated step-size selection. ACVI is noisy and shows no clear convergence. In terms of feasibility, both P-OpEx (without constraint reformulation) and ACVI perform well,
although ACVI remains noisy. Because our methods do not enforce all constraints simultaneously, a residual feasibility gap persists
(Figure~\ref{fig_infeas_gap_ply1}) but decreases with the number of
iterations. Its behavior depends on the players' optimal decisions. In our
experiment, the matrix $A$ yields a strictly feasible optimal decision for
Player~2, whereas Player~1's optimum lies on the boundary. Consequently,
Player~2's feasibility violation decreases faster, while Player~1's ACVI trajectory exhibits substantial noise because the log-barrier penalty of this interior-point method is applied near a boundary solution. P-OpEx with constraint reformulation (both with the same and different samples $\bar \omega_k^i$ and $\omega_k^i$) converges slowly.



\section{Conclusion}

This paper presents modified stochastic Korpelevich and Popov methods with randomized feasibility updates for solving monotone stochastic variational inequalities with a large (possibly infinite) number of constraints, expressed as the intersection of convex functional level sets. We establish the convergence rate of $O(1/\sqrt{T})$ for the modified dual-gap function, matching that of the existing primal–dual and ADMM-based approaches. Numerical results for a zero-sum matrix game further indicate favorable performance of our methods relative to the existing methods.


\section*{Acknowledgments}
We thank the reviewers for their constructive feedback.


\bibliographystyle{siamplain}
\bibliography{references}


\appendix

\section{Issue with constraint qualification assumption in~\cite{tang2023optimal}}\label{app_flaw}

Work \cite{tang2023optimal} considers the reformulation of the feasible set
\[
    X=\bigcap_{a\in\mathcal{A}}
    \left\{x \mid g_a(x)\leq 0\right\}
    =
    \left\{
        x \;\middle|\;
        \bar{\psi}(x)
        := \mathsf{E}_{a\sim\mathcal{A}}
        \left[g_a^+(x)\right]
        \leq 0
    \right\},
\]
which results in a degenerate inequality constraint. In fact, $\bar{\psi}(x)\geq 0$ for all $x$ by construction, so every point $x\in X$ must satisfy $\bar{\psi}(x)=0$. Therefore, the reformulated constraint has no strictly feasible (Slater) point. This degeneracy prevents the direct application of standard primal-dual convergence arguments that rely on strict feasibility to establish the existence and boundedness of Lagrange multipliers and the corresponding KKT guarantees. We note that~\cite[Theorem~2]{tang2023optimal} assumes Mangasarian-Fromovitz
constraint qualification (MFCQ)
holds at all limit points of the reformulated problem. For our case, the set $S := X \cap Y$. Hence, by~\cite[Remark~1]{tang2023optimal}, we require
\begin{equation}
    \min_{s\in\partial\bar{\psi}(x^*)} s^\top \nu > 0,
    \;\; \hbox{for some }
    \nu\in N_Y^*(x^*) , \tag{MFCQ} \label{MFCQ}
\end{equation}
at a feasible point $x^*\in S$ satisfying $\bar{\psi}(x^*)=0$, where $N_Y^*(x^*)$ is the the dual cone of the normal cone
$N_Y(x^*)$. Note that \cite[Remark 1]{tang2023optimal} has a typographical omission where $\nu \in N_Y(x^*)$ is considered instead of $N_Y^*(x^*)$. However, since
$\bar{\psi}(x)\geq 0$ for all $x$, the point $x^*$ is a global minimizer of $\bar{\psi}$, and therefore $0\in\partial\bar{\psi}(x^*)$. Hence, the condition in~\eqref{MFCQ} fails to hold. Consequently, the strict inequality in~\eqref{MFCQ} cannot hold either in the setting of \cite{tang2023optimal} or in our setting under this reformulation. Therefore, the standard MFCQ is incompatible with the exact positive-part expectation reformulation, despite being assumed in~\cite{tang2023optimal}.


\section{P-OpEx~\cite[Algorithm~2]{boob2026first} with reformulated constraints}\label{app_popex}

We implement P-OpEx~\cite[Algorithm~2]{boob2026first} using the single expectation
constraint $\{x \mid \mathsf{E}_{a\sim\A}\!\left[g_a^+(x)\right]\leq 0\}$, with one associated Lagrange multiplier $\lambda\in\R_+$:
\begin{align*}
\lambda_{k+1}
&=
\Pi_{\lambda\geq 0}
\left[
\lambda_k+
\frac{1}{\tau_k}
\mathsf{E}_{a\sim\A}\!\left[g_a^+(x_k)\right]
\right],\\
x_{k+1}
&=
\Pi_Y
\!\left[
x_k-\frac{1}{\eta_k}
\left(
(1+\theta_k)\widehat F(x_k,\xi_k)
-\theta_k\widehat F(x_{k-1},\xi_{k-1})
\!+\! \lambda_{k+1}
\partial_x
\mathsf{E}_{a\sim\A}\!\left[g_a^+(x_k) \right]
\right)
\!\right],
\end{align*}
where $\tau_k$ and $\eta_k$ are the inverse step sizes and $\theta_k$ is the
extrapolation parameter. $g_a^+(x)$ is convex for all $a\in\A$ and we assume that there exists an open set $U\supseteq Y$ such that, for
every $x\in U$, the mapping $a\mapsto g_a^+(x)$ is measurable and $\mathsf{E}_{a\sim P}\!\left[g_a^+(x)\right]<\infty$. Under these conditions, \cite{bertsekas1972stochastic} allows $\partial_x
\mathsf{E}_{a\sim\A}\!\left[g_a^+(x_k)\right]
=
\mathsf{E}_{a\sim\A}\!\left[\partial_x g_a^+(x_k)\right]$. Hence, expectations can be estimated using stochastic
samples as shown in Algorithm~\ref{alg:single-expectation-popex}.
%
%
\begin{algorithm}[H]
\caption{Partial Operator Extrapolation for the expectation
reformulation}
\label{alg:single-expectation-popex}
\begin{algorithmic}[1]
\REQUIRE $x_0\in Y$, $\lambda_0=0$.
\STATE Set $x_{-1}=x_0$, sample $\xi_{-1} \in \Xi$.
\FOR{$k=0,1,\ldots,T-1$}
    \FOR{$i=1,\ldots,N_{k}$}
        \STATE Independently sample $\bar \omega_k^i, \omega_k^i \sim \A$ and $\xi_k \in \Xi$
        \STATE Set $h_k^i = g_{\bar \omega_k^i}^+ (x_k)$ and $s_k^i = \partial_x g_{\omega_k^i}^+(x_k)$
        \STATE Update: $\lambda_{k+1} = \Pi_{\lambda \geq 0} [\lambda_k + \frac{1}{\tau_k} h_k^i]$
        \STATE Update: $x_{k+1} = \Pi_Y [x_k - \frac{1}{\eta_k} [(1+\theta_k) \widehat F(x_k, \xi_k) - \theta_k \widehat F(x_{k-1}, \xi_{k-1}) + \lambda_{k+1} s_k^i]]$
    \ENDFOR
\ENDFOR
\STATE \textbf{return} $\bar x_T = \frac{1}{T} \sum_{t=0}^{T-1}x_{t+1}$
\end{algorithmic}
\end{algorithm}
%
%
%
%
The method uses $N_k$ samples in each iteration. The samples $\bar\omega_k^i,\omega_k^i\in\A$ may be identical or distinct.


\end{document}